\documentclass[12pt]{article}
\usepackage{amsmath,amssymb,amsthm}
\usepackage{mathtools}
\usepackage{graphics,epsfig}
\usepackage{hyperref}
\usepackage{natbib}
\usepackage{color}
\usepackage{graphicx}
\usepackage{caption}
\usepackage{subcaption}
\usepackage{float}
\usepackage{mathrsfs}
\usepackage{multirow}
\usepackage{comment}
\usepackage{setspace}
\usepackage{bbm}
\usepackage{bm}
\usepackage{amsfonts}
\usepackage{xcolor}
\usepackage{pslatex}
\usepackage{setspace}
\usepackage{bbm}

\def \ra {\rightarrow}

\def \E {\mathbb{E}}

\def \a {\alpha}
\def \be {\beta}

\newtheorem{definition}{\bf Definition}

	\newtheorem{remark}{\bf Remark}

	\newtheorem{theorem}{\bf Theorem}
	\newtheorem{thm}[theorem]{\bf Theorem}
	\newtheorem{prop}[theorem]{\bf Proposition}
	\newtheorem{lem}[theorem]{\bf Lemma}
	\newtheorem{cor}[theorem]{\bf Corollary}

	\newtheorem{as}[theorem]{\bf Assumption}
	\newcommand{\sgn}{\operatorname{sgn}}


\begin{document}

  \title{\bfseries An optimal level of Stubbornness to win a soccer match}

  \author{
Paramahansa Pramanik\footnote{Corresponding author, {\small\texttt{ppramanik@southalabama.edu}}}\; \footnote{Department of Mathematics and Statistics,  University of South Alabama, Mobile, AL, 36688,
United States.}
}

\date{\today}
\maketitle

\begin{abstract}
This study conceptualizes stubbornness as an optimal feedback Nash equilibrium within a dynamic setting. To assess a soccer player's performance, we analyze a payoff function that incorporates key factors such as injury risk, assist rate, passing accuracy, and dribbling ability. The evolution of goal-related dynamics is represented through a backward parabolic partial stochastic differential equation (BPPSDE), chosen for its theoretical connection to the Feynman-Kac formula, which links stochastic differential equations (SDEs) to partial differential equations (PDEs). This relationship allows stochastic problems to be reformulated as PDEs, facilitating both analytical and numerical solutions for complex systems. We construct a stochastic Lagrangian and utilize a path integral control framework to derive an optimal measure of stubbornness. Furthermore, we introduce a modified Ornstein-Uhlenbeck BPPSDE to obtain an explicit solution for a player's optimal level of stubbornness.
\end{abstract}

\subparagraph{Key words:} Stochastic differential games; BPPSDE;  Feedback Nash equilibrium; sports analytics.

\section{Introduction.}
This study determines the optimal degree of stubbornness for a soccer player navigating goal-scoring dynamics, modeled through a backward parabolic partial stochastic differential equation. A Feynman-type path integral control method is employed to derive this measure. Stubbornness is a crucial trait in soccer, particularly for goal-scoring, as it fosters perseverance and resilience. Several key aspects underscore its benefits: first, tenacious players continue their efforts despite strong defenses or repeated failures, increasing their chances of creating scoring opportunities. Second, obstacles such as defensive pressure, harsh weather, or fatigue may hinder success, but persistent players remain focused until they achieve their goal \citep{pramanik2023scoring,pramanik2024motivation}. Third, maintaining a determined mindset strengthens a player’s confidence, reinforcing their belief in eventually overcoming the goalkeeper. Fourth, when paired with adaptability, stubbornness allows players to learn from missed chances, refining their approach rather than becoming discouraged, ultimately leading to greater success \citep{pramanik2020motivation}. Finally, relentless determination can inspire teammates, elevating overall team morale and improving the likelihood of breakthroughs. In this study, we consider how a player’s level of stubbornness is influenced by a mean field framework. Additionally, we derive an explicit solution for stubbornness as a feedback Nash equilibrium, incorporating key performance factors such as injury risk, assist rate, passing accuracy, and dribbling ability.

Soccer enjoys widespread popularity globally, thanks to its straightforward rules. Among the most significant tournaments in the sport are the World Cup, the Euro Cup, and the Copa America, with the World Cup being the most renowned. \cite{santos2014} emphasized that FIFA has been concerned about teams becoming overly defensive in their goal-scoring strategies since the early 1990s, which has led to a reduction in the total number of goals over time and a subsequent dip in the sport’s appeal. For example, during the 2006 men's World Cup, both Italy and Spain conceded only two goals throughout their seven matches. A USA Today article from March 17, 1994, highlighted FIFA's goal to promote more attacking and high-scoring games. This objective has prompted some teams to alter their approach, as seen in the 2014 men's World Cup semifinal, where Germany scored five goals against Brazil in the first half \citep{pramanik2024dependence,pramanik2024parametric}. In the 2018 men's World Cup, France adopted an offensive strategy despite lacking star players and having a younger team, allowing them to play freely and ultimately win the tournament. On the other hand, teams with well-known soccer stars sometimes underperformed because their predictable playstyles, honed over long careers, were easily countered by opponents. This highlights how mental stress can significantly impact game outcomes, influencing the financial stakes for both winners and losers \citep{pramanik2021optimization,pramanik2021optimal,pramanik2023scoring}.

Soccer is a game rich in strategic complexity, providing numerous avenues for analyzing optimal tactics \citep{dobson2010}. In a match, two teams face off for ninety minutes, each composed of a goalkeeper and ten outfield players. Coaches have the flexibility to organize their players in any formation and can alter their team's playing style at any point during the game \citep{dobson2010}. The primary goal for each team is to score while simultaneously preventing the opposition from doing so \citep{pramanik2024parametric}. The strategies employed by both teams are highly interdependent, as the approach taken by one team directly affects the other’s chances of scoring or conceding goals \citep{dobson2010}.

Palomino (1998) made a notable impact on the development of dynamic game-theoretic models aimed at optimizing strategic decisions in sports \citep{palomino1998}. This research examined how two competing soccer teams continuously chose between defensive and offensive formations. It was observed that a team leading in the score tends to adopt a defensive stance when the match is tied or they are behind, particularly in the second half. In a similar vein, \cite{banerjee2007} explored the strategic shifts in the National Hockey League following a change in the points scoring system for games tied at the end of regulation time \citep{dobson2010}. The analysis in these studies, however, is constrained by certain assumptions that are either unrealistic or contentious, limiting the scope for a broader dynamic strategic evaluation \citep{pramanik2020optimization,pramanik2024optimization}. For example, one study assumes—without verification—that adopting an offensive strategy boosts the goal-scoring rate of the attacking team more significantly than that of the defending team. This assumption is tested by examining the 2014 men's soccer World Cup semifinal between Germany and Brazil, where Germany scored five goals in just 29 minutes, including four within the first six minutes, ultimately leading 7-0 by the second half. However, despite Germany's aggressive approach early on, they only managed to score two additional goals in the remaining hour of the game. Another study challenges the previous assumption by introducing a parameter that accounts for a comparative advantage in either offense or defense \citep{polansky2021motif,pramanik2024estimation}. Both articles under review focus on the strategic decisions available to soccer teams, yet they primarily discuss the binary options of attack and defense, neglecting a vital secondary dimension: the choice between a violent and non-violent style of play \citep{dobson2010,hua2019assessing}. Teams that play violently may commit fouls and risk red cards, while other aggressive acts may aim to disrupt or sabotage the opposing team \citep{dobson2010}. The study by \cite{dobson2010} assumes that these strategic choices are discrete, forcing teams to select between defensive and attacking formations, as well as violent and non-violent styles of play. These choices significantly affect the chances of scoring or conceding goals and the probability of receiving red cards. Through numerical simulations, Dobson shows that optimal strategic decisions depend on various factors, including the current score difference and the time remaining in the game \citep{dobson2010}. In this paper, we treat stubbornness as a continuous control variable, with its maximum value representing an attacking approach \citep{pramanik2022lock,pramanik2022stochastic}.

The structure of the paper is as follows: Section 2 discusses the properties of the payoff function and the BPPSDE. Section 3 presents an explicit solution for optimal stubbornness using the Wick-rotated Schr\"odinger-type equation. Finally, Section 4 offers a brief conclusion.

\section{Background Framework.}

In this section, we will be discussing the construction of a forward-looking stochastic goal dynamics along with a conditional expected dynamic objective function.

\subsection{Probabilistic Construction.}

Let $t>0$ be a fixed finite time. Define $W_s$ is a one-dimensional Wiener process appeared in a soccer game at time $s$ for all $s\in[0,t]$ defined on a complete probability space $(\Omega, \mathcal F,\mathbb P)$ with sample space $\Omega$, $\sigma$-field and probability measure $\mathbb P$.

\begin{definition}\label{d0}
Consider $\left\{\mathcal F_t\right\}_{t\in I}$	as a family of sub-$\sigma$-fields of $\mathcal F$, such that $I$ represents an ordered indexed set with the condition $\mathcal F_s\subset\mathcal F_t$ for every $s<t$, and $s,t\in I$. This family $\left\{\mathcal F_t\right\}_{t\in I}$ is called a filtration of the above process.
\end{definition}

If we simply talk of a goal dynamics $\{x\}_{t\in I}$ or simply $x(t)$, then this implies the choice of filtration corresponding to a soccer match is
\[
\mathcal F_t:=\mathcal F_t^x:=\sigma\left\{x(s)|s\leq t, s,t\in I\right\},
\] 
which is called as canonical or natural filtration of ${x}_{t\in I}$. In our case, the Wiener process is associated with the canonical filtration,
\[
\mathcal F_t^W:=\sigma\left\{W_s|0\leq s\leq t\right\},\ \ t\in[0,\infty).
\]
In this paper, $x(s)\in \mathcal X$ denotes the goal dynamics at time $s$ for an individual player, with values in $\mathbb R$. The goalkeepers are excluded from this analysis, as they typically remain at their goal posts and generally do not move to the opposite side to score. By goal dynamics, we refer to the probability of scoring a goal at a specific continuous time point $s$. A higher probability of scoring does not guarantee that a goal has actually been scored \citep{pramanik2024semicooperation,pramanik2023path1,pramanik2023path}. For instance, this includes situations where a striker misses the goal, a shot is deflected by the opposing goalkeeper or players, or an offside call is made, among other possibilities. Let the goal dynamics be represented by a stochastic differential equation (SDE)
\begin{equation}\label{0}
dx(s)=\mu\left[s,x(s),u(s),\sigma_2\left[s,x(s),u(s)\right]\right]ds+\sigma\left[s,x(s),u(s)\right]dW_s,
\end{equation}
where $\mu:[0,t]\times\mathcal X\times\mathcal U\times\mathbb R\mapsto \mathbb R$, $\sigma_2:[0,t]\times\mathcal X\times\mathcal U\mapsto \mathbb R$ and $\sigma:[0,t]\times\mathcal X\times\mathcal U\mapsto \mathbb R$ are the drift and diffusion coefficients, respectively, $x(0)=x_0$, and the control variable $u\in\mathcal U$ takes values in $\mathbb R$. In this context, the control variable $u$ represents the level of stubbornness of a player. A highly stubborn player makes decisions independently, disregarding strategies set in the dressing room or by the coach, leading to a high value of $u$. Conversely, a less stubborn player adheres to team strategies, resulting in a low value of $u$. We can assume, without loss of generality, that $u \in [0,1]$, where $0$ indicates no stubbornness \citep{pramanik2016tail,pramanik2021effects}. Furthermore, define
\begin{equation}\label{1}
\sigma\left[s,x(s),u(s)\right]:=\sigma_1\left[s,x(s),u(s)\right]-\sigma_2\left[s,x(s),u(s)\right],
\end{equation}
where $\sigma_1$ represents a stochastic network of passing the ball from a single player to others. In this context, we account for the passing of 19 players, as a missed pass might result in the ball being intercepted by an opposing player, and $\sigma_2$ arises due to environmental and strategic uncertainties \citep{pramanik2024motivation}. The negative sign in Equation \eqref{1} indicates that $\sigma_2$ negatively impacts $\sigma_1$, as environmental and strategic uncertainties contribute to reducing the size of the stochastic ball-passing network. The SDE of the form of Equation \eqref{1} is called as one backward parabolic partial stochastic differential equation (BPPSDE).

\subsection{Existence and Uniqueness of a solution of the BPPSDE.}

In this section, we will demonstrate the existence and uniqueness of the BPPSDE described in Equation \eqref{1}. In the following lemma, we adopt a similar approach to the Feynman-Kac formula. However, since we assume the existence of a unique solution to a stochastic parabolic partial differential equation, we cannot label it as the Feynman-Kac formula, as that specifically relates a deterministic linear parabolic partial differential equation to path integrals.

\begin{lem}
	Let
	\begin{align}\label{2.4}
	&\biggr[-\mathcal{V}[s,x(s),u(s)]\ \phi^*[s,x(s),u(s)]\notag\\
	&\hspace{.5cm}+\Theta [s,x(s),u(s)]+\frac{\partial \phi^*[s,x(s),u(s)]}{\partial s}\notag\\
	&\hspace{1cm}+\frac{\partial \phi^*[s,x(s),u(s)]}{\partial x}\mu\left\{s,x(s),u(s),\sigma_2 [s,x(s),u(s)]\right\}\notag\\
	&\hspace{1.5cm}+\frac{\partial^2 \phi^*[s,x(s),u(s)]}{2\partial x^2}\left\{\sigma_1 [s,x(s),u(s)]-\sigma_2 [s,x(s),u(s)]\right\}^2\biggr]ds\notag\\
	&\hspace{2cm}-\frac{\partial \phi^*[s,x(s),u(s)]}{\partial x}\sigma_2 [s,x(s),u(s)]\ dW_s=0,
	\end{align}
	be a Cauchy's problem of a BPPSDE where $[s,x(s),u(s)]\in[0,t]\times\mathbb{R}\times[0,1]$ subject to the terminal condition $\phi^*[t,x(t),u(t)]=\mathcal{T}[t,x(t)]$,  $\Theta,\mathcal{V},\mu,\sigma_1$, and $ \sigma_2$ are known functions and $\phi^*[s,x(s),u(s)]:[0,t]\times\mathbb{R}\times[0,1]\ra\mathbb{R}$ is an unknown function. Then, assuming a unique solution exists, it  can be written as the conditional expectation
	
	\begin{align}
	&\phi^*[s,x(s),u(s)]\notag\\&=\E\biggr\{\mathcal{T}[t,x(t)]\ \exp\left\{-\int_s^t\ \mathcal{V}[\kappa,x(\kappa),u(\kappa)]\ d\kappa\right\}\notag\\&\hspace{.5cm} + \int_s^t \Theta [s_1,x(s_1),u(s_1)]\ \exp\ \left\{-\int_s^{s_1}\ \mathcal{V}[\kappa,x(\kappa),u(\kappa)]\ d\kappa\right\}\ ds_1\ |X(s) = x(s),U(s)=u(s)\biggr\}\notag
	\end{align}
	for all $s_1\in[0,s]$ such that $ X$ be an It$\hat o$'s process of the form represented by SDE \eqref{0}.
	\label{lem1}
\end{lem}

\begin{proof}
	Assume that $\phi^*[s,x(s),u(s)]$ is a solution to the stochastic differential equation given in Equation (\ref{2.4}). Define a stochastic process
	\begin{align}\label{2.5}
	\mathcal{Z}(s_1)&=\exp\left\{-\int_s^{s_1}\ \mathcal{V}[\kappa,x(\kappa),u(\kappa)]\ d\kappa\right\}\phi^*[s_1,x(s_1),u(s_1)]\notag\\&\hspace{1cm} + \int_s^{s_1} \exp\left\{-\int_s^l\ \mathcal{V}(\kappa,x(\kappa),u(\kappa))d\kappa\right\}\ \Theta[l,x(l),u(l)]\ dl,
	\end{align}
	Applying It\^o's Formula to the first term on the right hand side of Equation (\ref{2.5}) yields,
	\begin{align}\label{2.6}
	d\mathcal{Z}(s_1)&=\left\{d\exp\left\{-\int_s^{s_1}\ \mathcal{V}[\kappa,x(\kappa),u(\kappa)]\ d\kappa\right\}\right\}\phi^*[s_1,x(s_1),u(s_1)]\notag\\
	 &\hspace{1cm} +\exp\left\{-\int_s^{s_1}\ \mathcal{V}[\kappa,x(\kappa),u(\kappa)]\ d\kappa\right\} \left\{d\phi^*[s_1,x(s_1),u(s_1)]\right\}\notag\\
	  & \hspace{2cm}+\left\{d\ \exp\left\{-\int_s^{s_1}\ \mathcal{V}[\kappa,x(\kappa),u(\kappa)]\ d\kappa\right\}\right\} \left\{d\phi^*[s_1,x(s_1),u(s_1)]\right\}\notag\\ &\hspace{3cm}+d\left\{\int_s^{s_1} \exp\left\{-\int_s^l\ \mathcal{V}[\kappa,x(\kappa),u(\kappa)]\ d\kappa\right\}\ \Theta [l,x(l),u(l)]\ dl\right\}.
	\end{align}
	
	Since
	\begin{align}\label{2.60}
	& d\exp\left\{-\int_s^{s_1}\ \mathcal{V}[\kappa,x(\kappa),u(\kappa)]\ d\kappa\right\}\notag\\&=-\mathcal{V}[s_1,x(s_1),u(s_1)]\ \exp\left\{-\int_s^{s_1}\ \mathcal{V}(\kappa,x(\kappa),u(\kappa))d\kappa\right\}\ ds_1,
	\end{align}
	we have that,
	\begin{align}\label{2.7}
	&\left\{d\exp\left\{-\int_s^{s_1}\ \mathcal{V}[\kappa,x(\kappa),u(\kappa)]d\kappa\right\}\right\} \left\{d\phi^*[s_1,x(s_1),u(s_1)]\right\}\notag\\ &=\left\{-\exp\left\{-\int_s^{s_1}\ \mathcal{V}[\kappa,x(\kappa),u(\kappa)]\ d\kappa\right\}\mathcal{V}[s_1,x(s_1),u(s_1)]\ ds_1\right\}\left\{d\phi^*[s_1,x(s_1),u(s_1)]\right\}=0,
	\end{align}
	as $ds_1\ d\phi^*=0$ \citep{oksendal2013stochastic}. The fourth term of equation (\ref{2.6}) becomes
	\begin{align}\label{2.8}
	& d\left\{\int_s^{s_1} \exp\left\{-\int_s^l\ \mathcal{V}[\kappa,x(\kappa),u(\kappa)]\ d\kappa\right\}\ \Theta [l,x(l),u(l)]\ dl\right\}\notag\\&= \exp\left\{-\int_s^{s_1}\ \mathcal{V}(\kappa,x(\kappa),u(\kappa))d\kappa\right\}\ \Theta [s_1,x(s_1),u(s_1)]\ ds_1.
	\end{align}
	Equations (\ref{2.6}) - (\ref{2.8}) imply
	\begin{align}\label{2.9}
	d\mathcal{Z}(s_1)&=-\exp\left\{-\int_s^{s_1}\ \mathcal{V}[\kappa,x(\kappa),u(\kappa)]\ d\kappa\right\}\mathcal{V}[s_1,x(s_1),u(s_1)]\ \phi^*[s_1,x(s_1),u(s_1)]\ ds_1\notag\\& \hspace{1cm}+\exp\left\{-\int_s^{s_1}\ \mathcal{V}[\kappa,x(\kappa),u(\kappa)]\ d\kappa\right\} \left\{d\phi^*[s_1,x(s_1),u(s_1)]\right\}\notag\\ &\hspace{1.5cm}+\exp\left\{-\int_s^{s_1}\ \mathcal{V}[\kappa,x(\kappa),u(\kappa)]\ d\kappa\right\}\ \Theta [s_1,x(s_1),u(s_1)]\ ds_1\notag\\&=\exp\left\{-\int_s^{s_1}\ \mathcal{V}[\kappa,x(\kappa),u(\kappa)]\ d\kappa\right\}\biggr\{d\phi^*[s_1,x(s_1),u(s_1)]+\Theta [s_1,x(s_1),u(s_1)]\ ds_1\notag\\&\hspace{1cm}-\mathcal{V}[s_1,x(s_1),u(s_1)]\ \phi^*[s_1,x(s_1),u(s_1)]\ ds_1\biggr\}.
	\end{align}
	
	By Lemma (\ref{ito}) of the Appendix it is shown that
	\begin{align}\label{2.10}
	d\phi^*[s_1,x(s_1),u(s_1)]&=\biggr\{\frac{\partial \phi^*[s_1,x(s_1),u(s_1)]}{\partial s_1}+\frac{\partial \phi^*[s_1,x(s_1),u(s_1)]}{\partial x}\mu\left\{s_1,x(s_1),u(s_1),\sigma_2 [s_1,x(s_1),u(s_1)]\right\}\notag\\&\hspace{.5cm}+\frac{1}{2}\frac{\partial^2 \phi^*[s_1,x(s_1),u(s_1)]}{\partial x^2}\left\{\sigma_1 [s_1,x(s_1),u(s_1)]-\sigma_2 [s_1,x(s_1),u(s_1)]\right\}^2\biggr\}\ ds_1\notag\\& \hspace{1cm}+\frac{\partial \phi^*[s_1,x(s_1),u(s_1)]}{\partial x}\sigma_1 [s_1,x(s_1),u(s_1)]\ dW_{s_1}\notag\\&\hspace{1.5cm}-\frac{\partial \phi^*[s_1,x(s_1),u(s_1)]}{\partial x}\sigma_2 [s_1,x(s_1),u(s_1)]\ dW_{s_1}.
	\end{align}
	Equations (\ref{2.10}) into (\ref{2.9}) imply,
	\begin{align}\label{2.11}
	& d\mathcal{Z}(s_1)\notag\\&=\exp\left\{-\int_s^{s_1}\ \mathcal{V}[\kappa,x(\kappa),u(\kappa)]\ d\kappa\right\}\biggr\{\big\{-\mathcal{V}[s_1,x(s_1),u(s_1)]\ \phi^*[s_1,x(s_1),u(s_1)]\notag\\&\hspace{.5cm}+\Theta [s_1,x(s_1),u(s_1)]+\frac{\partial \phi^*[s_1,x(s_1),u(s_1)]}{\partial s_1}\notag\\&\hspace{1cm}+\frac{\partial \phi^*[s_1,x(s_1),u(s_1)]}{\partial x}\mu\left\{s_1,x(s_1),u(s_1),\sigma_2 [s_1,x(s_1),u(s_1)]\right\}\notag\\&\hspace{1.5cm}+\frac{1}{2}\frac{\partial^2 \phi^*[s_1,x(s_1),u(s_1)]}{\partial x^2}\left\{\sigma_1[s_1,x(s_1),u(s_1)]-\sigma_2 [s_1,x(s_1),u(s_1)]\right\}^2\big\}ds_1\notag\\&\hspace{2cm}-\frac{\partial \phi^*[s_1,x(s_1),u(s_1)]}{\partial x}\sigma_2[s_1,x(s_1),u(s_1)]\ dW_{s_1}\biggr\}\notag\\&\hspace{2cm}+\exp\left\{-\int_s^{s_1}\ \mathcal{V}[\kappa,x(\kappa),u(\kappa)]\ d\kappa\right\}\frac{\partial \phi^*[s_1,x(s_1),u(s_1)]}{\partial x}\sigma_1 [s_1,x(s_1),u(s_1)]\ dW_{s_1}.
	\end{align}
	Because $\phi^*[s,x(s),u(s)]$ is a solution of the stochastic differential equation given in Equation (\ref{2.4}). Equation (\ref{2.11}) becomes,
	\begin{align}\label{2.12}
	d\mathcal{Z}(s_1)&=\exp\left\{-\int_s^{s_1}\ \mathcal{V}[\kappa,x(\kappa),u(\kappa)]\ d\kappa\right\}\frac{\partial \phi^*[s_1,x(s_1),u(s_1)]}{\partial x}\ \sigma_1 [s_1,x(s_1),u(s_1)]\ dW_{s_1}
	\end{align}
	The integral form of Equation (\ref{2.12}) is,
	\begin{align}\label{2.13}
	&\mathcal{Z}(t)-\mathcal{Z}(s)\notag\\&=\int_s^t\ \exp\left\{-\int_s^{s_1}\ \mathcal{V}[\kappa,x(\kappa),u(\kappa)]\ d\kappa\right\}\frac{\partial \phi^*[s_1,x(s_1),u(s_1)]}{\partial x}\sigma_1 [s_1,x(s_1),u(s_1)]\ dW_{s_1}
	\end{align}
\end{proof}

From Lemma (\ref{lem1}) we can conclude that if the Cauchy problem corresponding to BSPDE has a unique solution then that can be written as a conditional Expectation, and 

\begin{align}\label{2.17}
d\phi^*[s,x(s),u(s)]&=\biggr\{\mathcal{V}[s,x(s),u(s)]\ \phi^*[s,x(s),u(s)]\notag\\&\hspace{.5cm}-\Theta [s,x(s),u(s)]\notag\\&\hspace{1cm}-\frac{\partial \phi^*[s,x(s),u(s)]}{\partial x}\mu\left\{s,x(s),u(s),\sigma_2 [s,x(s),u(s)]\right\}\notag\\&\hspace{1.5cm}-\frac{1}{2}\frac{\partial^2 \phi^*[s,x(s),u(s)]}{\partial x^2}\left\{\sigma_1 [s,x(s),u(s)]-\sigma_2 [s,x(s),u(s)]\right\}^2\biggr\}\ ds\notag\\&\hspace{2cm}+\frac{\partial \phi^*[s,x(s),u(s)]}{\partial x}\sigma_2 [s,x(s),u(s)]\ dW_s.
\end{align}

In Lemma (\ref{lem1}), we demonstrate the relationship between a BPPSDE with a unique solution and the path integral. Since our main focus is on establishing the existence of a unique solution for a dynamic profit function governed by a stochastic differential equation, in Proposition (\ref{prop1}), we identify the conditions under which this stochastic differential equation has a unique solution \cite{pramanik2020optimization,pramanik2024optimization}. To achieve this, we first define a complete filtration, which is a sequence of complete $\sigma$-algebras in the $\mathcal{F}_s$-measurable space. Then, utilizing Assumptions (\ref{as2})-(\ref{as5}), and finally, with Definitions (\ref{def1}) and (\ref{def2}), we obtain the unique solution for the SDE \citep{pramanik2024semicooperation}.

For sample space $\Omega$, $\sigma$-algebra $\mathcal{F}$, and probability measure $\mathbb P$, let $(\Omega,\mathcal{F},\{\mathcal{F}\}_{s\geq 0},\mathbb{P})$ be a complete probability space of a $k$-dimensional Wiener process $W=\{W_s:s\geq 0\}$ such that $\mathcal{F}_{s\geq 0}$ is the natural filtration generated by the Wiener process, augmented by all $\mathbb{P}$-null sets in $\mathcal{F}$ \citep{pramanik2023optimal}. Suppose, for $t>0$, $\wp$ is the $\sigma$-field corresponding to the predictable sets on $\{\Omega,(0,t)\}$ associated with the filtration $\mathcal{F}_{s\geq 0}$, $L^p$ is the functional space which is defined using a natural generalization of the $p$-norm for finite dimensional vector spaces and $\tilde\partial_x^\a$ is the vector of weak derivatives with respect to the vector $x$ with $\tilde \a^{\text{th}}$ order, where $\tilde{\a} =(\a_1,\a_2,...,\a_\Im)'$ \citep{pramanik2021consensus,pramanik2024bayes}.

\begin{definition}
	(Sobolev space) For $\Im\geq1$ suppose $\Omega$ is an open set in $\mathbb{R}^\Im$. Let $p\geq1$ and $\jmath\in \mathbb{N}$. The Sobolev space $\mathcal{S}^{\jmath,p}$ is defined as 
	\begin{equation}\label{sob}
	\{f_s\in L^p(\Omega\times (0,t)\times\mathbb{R}^\Im);\ \text{for all}\  |\a|\leq\jmath,\ \partial_x^\a\in L^p(\Omega)\},
	\end{equation}
	where $\tilde{\a} =(\a_1,\a_2,...,\a_\Im)'$, $|\tilde \a|=\sum_{k=1}^\Im \a_k$ and $\tilde\partial_x^\a=(\partial_{x_1}^\a, \partial_{x_2}^\a,\partial_{x_3}^\a,...,\partial_{x_\Im}^\a)'$ is the vector of weak derivatives.
	\label{def1}
\end{definition}

Suppose $\jmath$ is an integer and $\mathcal{G}^\jmath=\mathcal{G}^\jmath(\mathbb{R}^\Im)$ be the Sobolev space $\mathcal{S}^{\jmath,p}(\mathbb{R}^\Im)$ on $\mathbb{R}^\Im$. Du et al. (2010) \cite{du2010} give conditions for finding a solution to the Equation (\ref{2.17}),

\begin{definition}
	The function pair $\left\{\phi^*[s,x(s),u(s)],\sigma_2[s,x(s),u(s)]\right\}\ \text{that maps}\ \Omega\times[0,t]\times\mathbb{R}$ to $\mathbb{R}^2$ is a weak solution to the parabolic partial stochastic differential Equation  (\ref{2.17}) if $\phi^*[s,x(s),u(s)]\in L^2[\Omega\times(0,t),\wp,\mathbb{G}^1]$ and $\sigma_2[s,x(s),u(s)]\in L^2[\Omega\times(0,t),\wp,\mathbb{G}^0]$, such that for each $\tau\in \mathbb{G}^1$ and each $(\omega,s)\in\Omega\times[0,t]$,
	
	\medskip
	
	\begin{align}\label{2.18}
	&\int_{\mathbb{R}}\ \phi^*[s,x(s),u(s)]\tau[x(s)]\ dx(s)\notag\\&=\int_{\mathbb{R}}\ \mathcal{T}[t,x(t)] \tau[x(t)]\ dx(t)+\int_{0}^t\int_{\mathbb{R}}\biggr\{-\mathcal{V}[s,x(s),u(s)]\ \phi^*[s,x(s),u(s)]+\Theta [s,x(s),u(s)]\notag\\&\hspace{.5cm}+\frac{\partial \phi^*[s,x(s),u(s)]}{\partial x}\mu\left\{s,x(s),u(s),\sigma_2 [s,x(s),u(s)]\right\}\notag\\&\hspace{1cm}+\frac{1}{2}\frac{\partial^2 \phi^*[s,x(s),u(s)]}{\partial x^2}\left\{\sigma_1 [s,x(s),u(s)]-\sigma_2 [s,x(s),u(s)]\right\}^2  \biggr\}\tau[x(s)]\ dx(s)\ ds\notag\\&\hspace{1.5cm}+\int_{0}^t\int_{\mathbb{R}}\ \frac{\partial \phi^*[s,x(s),u(s)]}{\partial x}\sigma_2 [s,x(s),u(s)] \ dx(s)\ dW_s
	\end{align}
	\label{def2}
\end{definition}

\begin{as}
	The functions $\phi^*,\mathcal{V},\Theta,\mu,\sigma_1$ and $\sigma_2$ are $\wp\times\mathbb{B}(\mathbb{R})$-measurable real valued functions where the terminal value $\mathcal{T}[t,x(t)]$ is a $\mathcal{F}_t\times\mathbb{B}(\mathbb{R})$ measurable real valued function, where $\mathbb{B}(\mathbb{R})$ is the Borel measurable set on $\mathbb R$.
	\label{as2}
\end{as}

\medskip

\begin{as}
	(Super parabolicity) For any two given constants $C\in(1,\infty)$ and $c\in(0,1)$ $c+\sigma_1^2(\omega,s,x)+\sigma_2^2(\omega,s,x)\leq2\sigma_1(\omega,s,x)\sigma_2(\omega,s,x)\leq C$  for all $(\omega,s,x)\in\Omega\times[0,t]\times\mathbb{R}^1$.
	\label{as3}
\end{as}

\begin{as}
	Suppose $\mathcal{H}_1$ and $\mathcal{H}_2$ are two separable Hilbert spaces such that $\mathcal{H}_1$ is densely embedded in $\mathcal{H}_2$. We also assume that both $\mathcal{H}_1$ and $\mathcal{H}_2$ have the same dual space $\mathcal{H_2}$ such that $\mathcal{H}_1\subset \mathcal{H}_2\subset \mathcal{H_2}$ with the linear operators $\Xi_1(\omega,s):\mathcal{H}_1\ra\mathcal{H_2}$ and $\Xi_2(\omega,s):\mathcal{H}_2\ra\mathcal{H_2}$ given by
	\begin{align}
	\Xi_1(\omega,s)\phi^*[s,x(s),u(s)]\ &= \frac{\partial \phi^*[s,x(s),u(s)]}{\partial x}\mu\left\{s,x(s),u(s),\sigma_2 [s,x(s),u(s)]\right\},\notag 
	\end{align}
	and
	\begin{align}
	\Xi_2(\omega,s)\sigma_i[s,x(s),u(s)]\ &= \frac{1}{2}\frac{\partial^2 \phi^*[s,x(s),u(s)]}{\partial x^2}\left\{\sigma_1 [s,x(s),u(s)]-\sigma_2 [s,x(s),u(s)]\right\}^2 \notag 
	\end{align}
	for $i=1,2$.
	We further assume $\mathcal{T}[t,x(t)]$ and $\Theta[s,x(s),u(s)]$ take values in $\mathcal{H}_1$ and $\mathcal{H_2}$.
	\label{as4} 
\end{as}

Denote the norms of the spaces $\mathcal{H}_1, \mathcal{H}_2$ and $\mathcal{H_2}$ as $||.||_{\mathcal{H}_1}, ||.||_{\mathcal{H}_2}$ and $||.||_{\mathcal{H_2}}$, respectively. Furthermore, in order to define inner products and duality products we use $(.\ ,.)$ and $\langle.\ ,.\rangle$, respectively. Finally, define 
\begin{align}
\hat{\sigma}[s,x(s),u(s)]&:=\left\{\frac{\partial}{\partial x} \phi^*[s,x(s),u(s)]\right\}\sigma_2 [s,x(s),u(s)]\in \mathcal{H}_2.
\end{align}

\begin{as}
	(coercivity condition)	For any two constants $\zeta_1>0$ and $\zeta_2>0$ there exists
	\begin{align}\label{2.19}
	2\langle x,\Xi_1x\rangle+||\Xi_2^*x||^2_{\mathcal{H}_2^*}&\leq-\zeta_1||x||^2_{\mathcal{H}_1}+\zeta_2||x||^2_{\mathcal{H}_2},\notag\\||\Xi_1 x||_{\mathcal{H}_2^*}&\leq\Xi_2||x||_{\mathcal{H}_1},
	\end{align}
	where $(\omega,s)\in\Omega\times[0,t]$ and $\Xi_2^*:\mathcal{H}_1\ra\mathcal{H}_2^*$ is the adjoint operator of $\Xi_2$.
	\label{as5}
\end{as}

\begin{prop}
	Suppose all the previous assumptions and definitions hold. Furthermore let us consider $\Theta\in L^2(\Omega\times(0,t),\wp,\mathcal{H_2})$ and $\mathcal{T}\in L^2(\Omega,\mathcal{F}_t,\mathcal{H}_2)$. Then Equation (\ref{2.17}) has a unique solution $(\phi^*,\hat \sigma)\in L^2(\Omega\times(0,t),\wp,\mathcal{H}_1\times\mathcal{H}_2^*)$ such that $\phi*\in \zeta([0,t],\mathcal{H}_2)$ (almost sure), and moreover,
	\begin{align}
	&\E\sup_{s\leq t}||\phi^*(s,x(s),u(s))||_{\mathcal{H}_2}^2+\E\int_0^t\ \big(||\phi^*[s,x(s),u(s)]||_{\mathcal{H}_1}+||\hat\sigma[s,x(s),u(s)]||_{\mathcal{H}_2^*}\big)\ ds\notag\\&\hspace{2cm}\leq\zeta\biggr( \E \int_0^t ||\Theta[s,x(s),u(s)]||_{\mathcal{H_2}}\ ds+\E||\mathcal{T}[t,x(t),u(t)]||_{\mathcal{H}_2^2}\biggr)\notag
	\end{align}
	where $\zeta=\zeta(\zeta_1,\zeta_2,t)$ is a constant. 
	\label{prop1}
\end{prop}

\begin{proof}
	To prove the above theorem we will take two steps. In the first step, we will assume that a solution to equation (\ref{2.17}) exists and demonstrate that this solution is unique. In the second step, we will establish the existence of the solution. This approach can be used in cancer research \citep{dasgupta2023frequent,hertweck2023clinicopathological,kakkat2023cardiovascular,khan2023myb,vikramdeo2023profiling,khan2024mp60,vikramdeo2024abstract}.
	
	Du et al. (2010) \cite{du2010} imply a process $(\phi^*,\hat{\sigma})$ is a $\mathcal{F}_s$-adapted $\mathcal{H}_1\times\mathcal{H}_2^*$ valued solution of equation (\ref{2.17}), if $\phi^*\in L^2(\Omega\times(0,t),\wp,\mathcal{H}_1)$ and  $\hat\sigma\in L^2(\Omega\times(0,t),\wp,\mathcal{H}_2^*)$ such that for each $\eta\in\mathcal{H}_1$ and $(\omega,s)\in\Omega\times[0,t]$ (almost everywhere),
	\begin{align}\label{2.20}
	(\eta,\phi^*(s,x(s),u(s)))&=\int_0^t\ [\langle\eta,\Xi_1\phi^*(s,x(s),u(s))\rangle+\Xi_2\ \hat{\sigma}(s,x(s),u(s))\notag\\&\hspace{1cm}+\Theta(s,x(s),u(s))]\ ds-\int_0^t\ \hat{\sigma}(s,x(s),u(s))\ dW_s.
	\end{align}
	Lemma 3.1 in \cite{du2010} implies $\phi^*(s,x(s),u(s))\in\zeta([0,t],\mathcal{H}_2)$ (almost sure).
	
	Let $\E\sup_{s\leq t}||\phi^*(s,x(s),u(s))||_{\mathcal{H}_2}^2$ be measurable. In other words, $\E\sup_{s\leq t}||\phi^*(s,x(s),u(s))||_{\mathcal{H}_2}^2<\infty$ and we define a stopping time
	\begin{align}\label{2.21}
	\aleph_k(\omega)&=\inf\{s;\ \sup_{s\leq t}||\phi^*(\omega,s,x(s),u(s))||_{\mathcal{H}_2}\geq k\}\wedge t.
	\end{align}
	Condition (\ref{2.21}) implies that $\aleph_k(\omega)\ra t$ almost surely as $k\ra\infty$. After implementing It\^o's formula Equation (\ref{2.17}) yields
	\begin{align}\label{2.22}
	\phi^*(t\wedge\aleph_k,x(t),u(t))&=-\int_0^{s\wedge\aleph_k}\big[-\mathcal{V}(s,x(s),u(s))\ \phi^*(s,x(s),u(s))\notag\\&\hspace{.5cm}+\Theta (s,x(s),u(s))\notag\\&\hspace{1cm}+\frac{\partial \phi^*(s,x(s),u(s))}{\partial x}\mu[s,x(s),u(s),\sigma_2 (s,x(s),u(s))]\notag\\&\hspace{1.5cm}+\frac{\partial^2 \phi^*(s,x(s),u(s))}{2\partial x^2}[\sigma_1 (s,x(s),u(s))-\sigma_2 (s,x(s),u(s))]^2\big]ds\notag\\&\hspace{2cm}+\int_o^{s\wedge\aleph_k}\frac{\partial \phi^*(s,x(s),u(s))}{\partial x}\sigma_2 (s,x(s),u(s))dW_s.
	\end{align}
	Moreover, Assumption (\ref{as4}) and lemma 3.1 in \cite{du2010} imply,
	\begin{align}\label{2.23}
	\phi^*(t\wedge\aleph_k,x(t),u(t))&=-\int_0^{s\wedge\aleph_k}\big[2\langle\phi^*(s,x(s),u(s)),\Xi_1\phi^*(s,x(s),u(s))\rangle\notag\\&\hspace{.5cm}+2(\Xi_2^*\phi^*(s,x(s),u(s)),\hat{\sigma}(s,x(s),u(s)))\notag\\&\hspace{1cm}+2\langle\phi^*(s,x(s),u(s)),\Theta(s,x(s),u(s))\rangle-||\hat{\sigma}(s,x(s),u(s))||_{\mathcal{H}_2^*}^2\big]\ ds\notag\\&\hspace{1.5cm}+\int_0^{s\wedge\aleph_k}2(\phi^*(s,x(s),u(s))
	,{\sigma_2}(s,x(s),u(s))\ dW_s,
	\end{align}
	and Assumption (\ref{as5}) yields
	\begin{align}\label{2.24}
	&\phi^*(t\wedge\aleph_k,x(t),u(t))\notag\\&=-\int_0^{s\wedge\aleph_k}\biggr[2\langle\phi^*(s,x(s),u(s)),\Xi_1\phi^*(s,x(s),u(s))\rangle\notag\\&\hspace{.5cm}+2(\Xi_2^*\phi^*(s,x(s),u(s)),\hat{\sigma}(s,x(s),u(s)))\notag\\&\hspace{1cm}+2\langle\phi^*(s,x(s),u(s)),\Theta(s,x(s),u(s))\rangle-||\hat{\sigma}(s,x(s),u(s))||_{\mathcal{H}_2^*}^2\biggr]\ ds\notag\\&\hspace{1.5cm}+\int_0^{s\wedge\aleph_k}2(\phi^*(s,x(s),u(s))
	,{\sigma_2}(s,x(s),u(s))\ dW_s\notag\\&\hspace{1cm}\leq\zeta(\zeta_2)\int_0^t\ \biggr(||\phi^*(s,x(s),u(s))||_{\mathcal{H}_1}^2+||\sigma_2(s,x(s),u(s))||_{\mathcal{H}_2^*}^2+||\Theta(s,x(s),u(s))||_{\mathcal{H}_1^*}^2\biggr)\ ds\notag\\&\hspace{1.5cm}+\int_0^{s\wedge\aleph_k}2(\phi^*(s,x(s),u(s))
	,{\sigma_2}(s,x(s),u(s))\ dW_s.
	\end{align}
	
	Burkholder-Davis-Gundy  (BDG) [see Proposition (\ref{prop14}) in appendix] implies
	\begin{align}\label{2.25}
	&\E\ \biggr|\sup_{s\leq\aleph_k}\ \int_0^t\ (\phi^*(s,x(s),u(s))
	,{\sigma_2}(s,x(s),u(s)))\ dW_s\biggr|\notag\\&\hspace{1cm}\leq\zeta\biggr[\E\ \int_0^{\aleph_k}\ ||\phi^*(s,x(s),u(s))||_{\mathcal{H}_2}^2\ ||\hat\sigma(s,x(s),u(s))||_{\mathcal{H}_2^*}^2 \biggr]^{\frac{1}{2}}\notag\\&\hspace{1cm}\leq\frac{1}{4}\ \E\ \sup_{s\leq\aleph_k}\ ||\phi^*(s,x(s),u(s))||_{\mathcal{H}_2}^2+\zeta\ \E\ \int_0^t\ ||\hat\sigma(s,x(s),u(s))||_{\mathcal{H}_2^*}^2\ ds.
	\end{align}
	
	Hence,
	\begin{align}\label{2.26}
	&\E\ \sup_{s\leq\aleph_k}\ ||\phi^*(s,x(s),u(s))||_{\mathcal{H}_2}^2\notag\\\hspace{1cm}& \leq\zeta(\zeta_2)\int_0^t\ \biggr(||\phi^*(s,x(s),u(s))||_{\mathcal{H}_1}^2+||\sigma_2(s,x(s),u(s))||_{\mathcal{H}_2^*}^2+||\Theta(s,x(s),u(s))||_{\mathcal{H}_1^*}^2\biggr)\ ds.
	\end{align}
	As our constant $\zeta$ does not depend on k, as $k\ra\infty$ we will get $\E\ \sup_{s\leq t}\ ||\phi^*(s,x(s),u(s))||_{\mathcal{H}_2}^2<\infty$.
	After using It\^o's formula one more time like in condition (\ref{2.24}) and using assumption (\ref{as5}) we get,
	\begin{align}\label{2.27}
	&||\phi^*(s,x(s),u(s))||_{\mathcal{H}_2}^2\notag\\&=\mathcal{T}(t,x(t),u(t))+\int_0^t\biggr[2\langle\phi^*(s,x(s),u(s)),\Xi_1\phi^*(s,x(s),u(s))\rangle\notag\\&\hspace{.5cm}+2(\Xi_2^*\phi^*(s,x(s),u(s)),\hat{\sigma}(s,x(s),u(s)))\notag\\&\hspace{1cm}+2\langle\phi^*(s,x(s),u(s)),\Theta(s,x(s),u(s))\rangle-||\hat{\sigma}(s,x(s),u(s))||_{\mathcal{H}_2^*}^2\biggr]\ ds\notag\\&\hspace{1.5cm}+\int_0^t 2(\phi^*(s,x(s),u(s))
	,{\sigma_2}(s,x(s),u(s))\ dW_s\notag\\&\leq \mathcal{T}(t,x(t),u(t))+ \int_0^t\biggr[2\langle\phi^*(s,x(s),u(s)),\Xi_1\phi^*(s,x(s),u(s))\rangle+(1+\epsilon)||\Xi_2^*\phi^*(s,x(s),u(s))||_{\mathcal{H}_2^*}^2\notag\\ &\hspace{.5cm}+\frac{1}{1+\epsilon}||\sigma_2(s,x(s),u(s))||_{\mathcal{H}_2^*}^2-||\sigma_2(s,x(s),u(s))||_{\mathcal{H}_2^*}^2\notag\\&\hspace{1cm}+\epsilon||\phi^*(s,x(s),u(s))||_{\mathcal{H}_1}^2+\frac{1}{\epsilon}||\Theta(s,x(s),u(s))||_{\mathcal{H_2}}^2\biggr]\ ds\notag\\&\hspace{1.5cm}+\int_0^t 2(\phi^*(s,x(s),u(s)),{\sigma_2}(s,x(s),u(s))\ dW_s\notag\\&\leq \mathcal{T}(t,x(t),u(t))+\int_0^t\biggr[-2\epsilon\langle\phi^*(s,x(s),u(s)),\Xi_1\phi^*(s,x(s),u(s))\rangle\notag\\&\hspace{.5cm}+(1+\epsilon)(-\zeta_1||\phi^*(s,x(s),u(s))||_{\mathcal{H}_1}^2+\zeta_2||\phi^*(s,x(s),u(s))||_{\mathcal{H}_2}^2)\notag\\&\hspace{1cm}-\frac{\epsilon}{1+\epsilon}||\sigma_2(s,x(s),u(s))||_{\mathcal{H}_2^*}^2+\epsilon||\phi^*(s,x(s),u(s))||_{\mathcal{H}_1}^2\notag\\&\hspace{1.5cm}+\frac{1}{\epsilon} ||\Theta(s,x(s),u(s))||_{\mathcal{H_2}}^2\biggr]\ ds+\int_0^t 2(\phi^*(s,x(s),u(s)),{\sigma_2}(s,x(s),u(s))\ dW_s\notag\\&\leq \mathcal{T}(t,x(t),u(t))+ \int_0^t\biggr[\{2\epsilon\zeta_2-\zeta_1(1+\epsilon)+\epsilon\}||\phi^*(s,x(s),u(s))||_{\mathcal{H}_1}^2+(1+\epsilon)\zeta_2||\phi^*(s,x(s),u(s))||_{\mathcal{H}_2}^2\notag\\&\hspace{.5cm}-\frac{\epsilon}{1+\epsilon}||\sigma_2(s,x(s),u(s))||_{\mathcal{H}_2^*}+\frac{1}{\epsilon}||\Theta(s,x(s),u(s))||_{\mathcal{H_2}}^2\notag\\&\hspace{1cm}+\int_0^t 2(\phi^*(s,x(s),u(s)),{\sigma_2}(s,x(s),u(s))\ dW_s.
	\end{align}
	Since $\epsilon\ra0$ then $2\epsilon\zeta_2-\zeta_1(1+\epsilon)+\epsilon\ra-\zeta_2<0$. Thus
	\begin{align}\label{2.28}
	&||\phi^*(s,x(s),u(s))||_{\mathcal{H}_2}^2+\int_0^t\big(||\phi^*(s,x(s),u(s))||_{\mathcal{H}_1}^2+||\sigma_2(s,x(s),u(s))||_{\mathcal{H}_2^*}\big)\ ds\notag\\&\hspace{.5cm}\leq \zeta(\zeta_1,\zeta_2)\biggr[\mathcal{T}(t,x(t),u(t))+\int_0^t\big( ||\phi^*(s,x(s),u(s))||_{\mathcal{H}_2}^2 +||\Theta(s,x(s),u(s))||_{\mathcal{H_2}}\big)\ ds\biggr]\notag\\&\hspace{1cm}+\int_0^t 2(\phi^*(s,x(s),u(s)),{\sigma_2}(s,x(s),u(s))\ dW_s.
	\end{align}
	
	As $\E\ \sup_{s\leq t}||\phi^*(s,x(s),u(s))||_{\mathcal{H}_2}^2<\infty$ repeating (\ref{2.27}) with $\int_0^t (\phi^*(s,x(s),u(s)),{\sigma_2}(s,x(s),u(s))\ dW_s$ is a uniform martingale. After taking expectations of both the sides in (\ref{2.28}) and by using Gronwall inequality [see corollary (\ref{cor16}) in appendix] we get,
	\begin{align}\label{2.29}
	&\sup_{s\leq t}\ \E\ ||\phi^*(s,x(s),u(s))||_{\mathcal{H}_2}^2+\E\ \int_0^t\big(||\phi^*(s,x(s),u(s))||_{\mathcal{H}_1}^2+||\sigma_2(s,x(s),u(s))||_{\mathcal{H}_2^*}\big)\ ds\notag\\&\hspace{.5cm}\leq\zeta\ e^{\zeta t}\ \biggr( \E\ \big[\mathcal{T}(t,x(t),u(t))\big]+ \E\ \int_0^t||\Theta(s,x(s),u(s))||_{\mathcal{H_2}}\ ds\biggr).
	\end{align}
	Equation (\ref{2.28}) and BDG inequality yield
	\begin{align}\label{2.30}
	&\E\ \sup_{s\leq t}||\phi^*(s,x(s),u(s))||_{\mathcal{H}_2}^2\notag\\&\hspace{.5cm}\leq\zeta(\zeta_1,\zeta_2,t)\biggr[\mathcal{T}(t,x(t),u(t))+ \int_0^t\ \big( ||\phi^*(s,x(s),u(s))||_{\mathcal{H}_2}^2 +||\Theta(s,x(s),u(s))||_{\mathcal{H_2}}\big)\ ds\biggr]\notag\\&\hspace{1cm}+\frac{1}{2}\E\ \sup_{s\leq t}||\phi^*(s,x(s),u(s))||_{\mathcal{H}_2}^2+\zeta\ \E\ \int_0^t\ ||\sigma_2(s,x(s),u(s))||_{\mathcal{H}_2^*}^2\ ds.
	\end{align}
	Inequality (\ref{2.30}) with (\ref{2.29}) gives the estimate of (\ref{2.25}).
	
	In this second part of the proof, we will demonstrate the existence of the solution. Assume that $\left\{b_i: i = 1, 2, 3, \dots\right\}$ is a complete orthogonal basis in $\mathcal{H}_2$, which also serves as an orthogonal basis in $\mathcal{H}_1$. We are employing the Galerkin approximation method here \citep{pramanik2024measuring}. Now, consider our system of BSDEs without the terminal condition in $\mathbb{R}^1$ with the basis function as follows
	\begin{align}\label{2.31}
	\phi_n^{*i}(s,x(s),u(s))&=\int_0^t\biggr[\langle b_i,\Xi_1(s,x(s),u(s))b_j\rangle\phi_n^{*j}(s,x(s),u(s))\notag\\&\hspace{.5cm}+\big(b_i,\Xi_2(s,x(s),u(s))\hat\sigma_n^j(s,x(s),u(s)) \big)\notag\\&\hspace{1cm}+\langle b_i,\Theta(s,x(s,u(s)))\rangle\biggr]\ ds-\int_0^t\sigma_{2n}^i(s,x(s),u(s))\ dW_s.
	\end{align}
	where $\phi_n^{*i}(s,x(s),u(s))$ and $\hat\sigma_n^j(s,x(s),u(s))$ are two unknown processes in $\mathbb{R}^1\times\mathbb{R}^1$. From \cite{du2010} we know that $\E\ [b_i,\mathcal{T}(t,x(t),u(t))]<\infty$ and $\E\ \int_0^t\langle b_i,\Theta(s,x(s),u(s))\rangle^2<\infty$. Then by \cite{du2010} there exists a unique solution of Equation \eqref{2.31}. Define this solution as $\phi_n^{*}(s,x(s),u(s)):=\sum_{i=1}^n\phi_n^{*i}(s,x(s),u(s))b_i$, and $\hat{\sigma}_n(s,x(s),u(s))=\sum_{i=1}^n\hat{\sigma}_n^i(s,x(s),u(s))b_i.$
	After applying It\^o's lemma and Burkholder-Davis-Gundy inequality of $||\phi_n^{*}(s,x(s),u(s))||_{\mathcal{H}_2}^2||$ like in the previous section we get,
	\begin{align}\label{2.32}
	&\E\ \sup_{s\leq t}\  ||\phi_n^*(s,x(s),u(s))||_{\mathcal{H}_2}^2+\E\ \int_0^t\big(||\phi_n^*(s,x(s),u(s))||_{\mathcal{H}_1}^2+||\sigma_{2n}(s,x(s),u(s))||_{\mathcal{H}_2^*}\big)\ ds\notag\\&\hspace{2cm}\leq\zeta(\zeta_1,\zeta_2)\biggr(\E\ [\mathcal{T}(t,x(t),u(t))]+\E\ \int_0^t||\Theta(s,x(s),u(s))||_{\mathcal{H_2}}^2\ ds\biggr).
	\end{align}
	From (\ref{2.32})  and \cite{du2010} we know $\exists\{n_1\}\in\{n\}$ where $\{n_1\}$ is the subsequence of $\{n\}$ and $(\phi^*,\hat{\sigma})\in L^2(\Omega\times(0,t),\wp,\mathcal{H}_1\times\mathcal{H_2})$ such that,
	\begin{align}\label{2.33}
	\phi_{n_1}^*\overset{w}\to \phi^*\in L^2(\Omega\times(0,t),\wp,\mathcal{H}_1),\notag\\\hat{\sigma}_{n_1}\overset{w}\to\hat{\sigma}\in L^2(\Omega\times(0,t),\wp,\mathcal{H_2}).
	\end{align}
	Suppose $\hbar:(\Omega,\mathcal{F})\ra\mathbb{R}^1$  and $\mho:[0,t]\ra\mathbb{R}^1$ where both of them are bounded measurable functions. Now the BSDE (\ref{2.31}) with $n\in\mathbb{N}$ and $b_i\in\{b_i\}$ where $i\leq n$ becomes,
	\begin{align}\label{2.34}
	&\E\ \int_0^t\ \hbar \mho(s)(b_i,\phi_{n_1}^*(s,x(s),u(s)))\ ds\notag\\&=\E\ \int_0^t\ \hbar \mho(s) \biggr\{\int_0^t\biggr[\langle b_i,\Xi_1(s,x(s),u(s))b_j\rangle\phi_n^{*j}(s,x(s),u(s))\notag\\&\hspace{.5cm}+\big(b_i,\Xi_2(s,x(s),u(s))\hat\sigma_n^j(s,x(s),u(s)) \big)\notag\\&\hspace{1cm}+\langle b_i,\Theta(s,x(s,u(s)))\rangle\biggr]\ ds-\int_0^t\sigma_{2n}^i(s,x(s),u(s))\ dW_s \biggr\}\ ds.
	\end{align}
	Furthermore,  \cite{du2010} imply,
	\begin{align}\label{2.35}
	\E\ \int_0^t\ \hbar \mho(s)(b_i,\phi_{n_1}^*(s,x(s),u(s)))\ ds\ra \E\ \int_0^t\ \hbar \mho(s)(b_i,\phi^*(s,x(s),u(s)))\ ds.
	\end{align}
	In view of the second condition of assumption \ref{as5} and estimate (\ref{2.32}), we get
	\begin{align}\label{2.36}
	\E\ \biggr|\hbar\langle b_i,\Xi_1\phi_{n_1}^*(s,x(s),u(s))\rangle \ ds \biggr|<\zeta<\infty.
	\end{align}
Since $\zeta$ is independent of $n_1$, for every $s\in[0,t]$
	\begin{align}\label{2.37}
	\E\ \int_0^t\ \hbar\langle b_i,\Xi_1\phi_{n_1}^*(s,x(s),u(s))\rangle \ ds\ra \E\ \int_0^t\ \hbar\langle b_i,\Xi_1\phi^*(s,x(s),u(s))\rangle \ ds.
	\end{align}
	Applying Fubini's theorem and Lebesgue dominated convergence theorem yield
	\begin{align}\label{2.38}
	&\E\ \int_0^t\ \hbar\mho(s)\int_{s_1}^t\ \langle b_i,\Xi_1\phi_{n_1}^*(s,x(s),u(s))\rangle \ ds_1ds\notag\\&= \int_0^t\ \mho(s)\ \E\ \int_{s_1}^t\ \hbar\langle b_i,\Xi_1\phi_{n_1}^*(s,x(s),u(s))\rangle \ ds_1ds\notag\\&\ra \int_0^t\ \mho(s)\ \E\ \int_{s_1}^t\ \hbar\langle b_i,\Xi_1\phi^*(s,x(s),u(s))\rangle \ ds_1ds.
	\end{align}
	Similarly,
	\begin{align}\label{2.39}
	&\E\ \int_0^t\ \hbar\mho(s)\int_{s_1}^t\ \big( b_i,\Xi_2\ \hat\sigma_{n_1}^*(s,x(s),u(s))\big) \ ds_1ds\notag\\&\hspace{1.5cm}\ra \E\ \int_0^t\ \hbar\mho(s)\int_{s_1}^t\ \big( b_i,\Xi_2\ \hat\sigma^*(s,x(s),u(s))\big) \ ds_1ds.
	\end{align}
	Assumptions (\ref{as5}) and Equation (\ref{2.32}) imply
	\begin{align}\label{2.40}
	\E\ \biggr|\hbar\int_0^t( b_i,\hat{\sigma}_{n_1}(s,x(s),u(s)))\ dW_s\biggr|<\zeta<\infty.
	\end{align}
	Since $(b_i,\hat{\sigma}_{n_1}(s,x(s),u(s)))\overset{w}\to(b_i,\hat{\sigma}(s,x(s),u(s))) $ in $L^2(0,t)$, previous result implies 
	\begin{align}\label{2.41}
	\int_0^t\ (b_i,\hat{\sigma}_{n_1}(s,x(s),u(s)))\ dW_s\overset{w}\to\ \int_0^t\ (b_i,\hat{\sigma}(s,x(s),u(s)))\ dW_s \in L^2(\Omega,\mathcal{F}_t,\mathbb{R}^1).
	\end{align}
Fubini's theorem and Lebesgue's dominated convergence theorem imply,
	\begin{align}\label{2.42}
	&\E\ \int_0^t\ \hbar\mho(s)\int_{s_1}^t\ \big( b_i, \hat\sigma_{n_1}^*(s,x(s),u(s))\big) \ dW_sds\notag\\&\hspace{1.5cm}\ra \E\ \int_0^t\ \hbar\mho(s)\int_{s_1}^t\ \big( b_i, \hat\sigma^*(s,x(s),u(s))\big) \ dW_sds.
	\end{align}
	Therefore, we find for $(\omega,s)\in\Omega\times[0,t]$ (a.e.),
	\begin{align}\label{2.43}
	(b_i,\phi^*(s,x(t),u(t))&=(b_i,\mathcal{T}(t,x(t),u(t)))+\int_0^t\biggr[\langle b_i,\Xi_1(s,x(s),u(s))b_j\rangle\phi^{*j}(s,x(s),u(s))\notag\\&\hspace{.5cm}+\big(b_i,\Xi_2(s,x(s),u(s))\hat\sigma^j(s,x(s),u(s)) \big)\notag\\&\hspace{1cm}+\langle b_i,\Theta(s,x(s,u(s)))\rangle\biggr]\ ds-\int_0^t\sigma_{2}^i(s,x(s),u(s))\ dW_s
	\end{align}
	Hence, we conclude that if we have the stochastic differential equation like BPPSDE, there exists a unique solution. This completes the proof.
	\end{proof}

In lemma \ref{lem1}, we demonstrate the connection between a unique solution to a BPPSDE and path integrals. Then, using assumptions \ref{as2}-\ref{as5} and proposition (\ref{prop1}), we establish the conditions under which a unique solution to the BPPSDE exists. Our goal is now to maximize the payoff function of a soccer player $\pi(s,x(s),u(s))$, where the stochastic goal dynamics take the form
\[
d x (s) =\mu[s,x(s),u(s),\sigma_2 (s,x(s),u(s))]\ ds + [\sigma_1 (s,x(s),u(s))- \sigma_2 (s,x(s),u(s))]\ dW_s.
\]
Let $x^*$ be the optimal output share of the firm corresponding to the optimal strategy $u^*(s,x^*(s))$. Following Du et al. (2013) \cite{du2013} define $\rho\in[0,1)$ such that a strategy function of spike variation becomes,
\begin{align}\label{2.44}
\tilde u(s,x(s)):=\begin{cases} u(s,x(s)) & \text{if $s\in[\rho,\rho+\epsilon)$}\\ u^*(s,x^*(s)) & \text{otherwise,} \end {cases}
\end{align}
where $\epsilon>0$ and $\epsilon\ra 0$ with any given strategy function $u(s,x(s))$. Condition (\ref{2.44}) implies $\tilde x(s)$ is the probability that a goal might be scored when strategy $\tilde u(s)$ has been used.

\begin{as}
	The stubbornness set $u(s,x(s)):[0,t]\times\mathcal{H}_2\ra\mathbb{R}^1$ is convex and the diffusion is independent of this set.\label{as7}
\end{as}

\begin{as}
	For each $(x,u)\in\mathcal{H}_2\times\mathcal{U}$ we have $\mu(.,x,u),\sigma_1(.,x,u),\sigma_2(.,x,u)$ and $\pi(.,x,u)$ all are predictable processes, where $\mathcal{U}$ is a non-empty Borel measurable subset of a metric space. We also assume for each $(s,x,u)\in(0,t)\times \mathcal{H}_2\times\mathcal{U}$; $\pi,\mu,\sigma_1$, and $\sigma_2$ are globally twice Frecht differentiable with respect to the market share $x$. The functions $\mu_x,\sigma_{1x},\sigma_{2x},\pi_{xx},\mu_{xx},\sigma_{1xx}$ are dominated by $\mathcal{M}$, $\sigma_{2xx}$ are continuous in $x(s)$ and dominated by $\mathcal{M}(1+||x||_{\mathcal{H}_2}+||u||_{\mathcal{H}_1})$. $\pi$ is dominated by $\mathcal{M}(1+||x||_{\mathcal{H}_2}^2+||u||_{\mathcal{H}_1}^2)$.
	\label{as8}
\end{as}

Suppose, $\mathcal{F}=\pi,\mu,\sigma_1,\sigma_2,\pi_x,\mu_x,\sigma_{1x},\sigma_{2x},\pi_{xx},\mu_{xx},\sigma_{1xx}$ such that

\begin{align}\label{2.45}
\tilde{\mathcal{F}}(s)&:=\mathcal{F}(s,\tilde x(s),\tilde u(s)),\notag\\\tilde{\mathcal{F}}^\Delta(s)&:=\mathcal{F}(s,\tilde x(s),\tilde u(s))-\tilde{\mathcal{F}}(s),\notag\\\tilde{\mathcal{F}}^\delta(s)&:=\tilde{\mathcal{F}}^\Delta(s). I_{[\rho,\rho+\epsilon]}(s),\notag\\\tilde{\mathcal{F}}_{xx}(s)&:=2\ \int_0^1\ \tau_1 \mathcal{F}_{xx} (s,\tau_1 x^*(s)+(1-\tau_1)\tilde x(s),\tilde u(s))\ d\tau_1.
\end{align}

\begin{thm}
	Let assumptions (\ref{as2}) through (\ref{as8}) hold. Define a penalized payoff function $\mathcal{P}:[0,t]\times\mathcal{H}_1\times[0,1]\ra\mathbb{R}^1$ as the form 
	\begin{align}\label{2.46}
	\mathcal{P}&:=\pi(s,x(s),u(s))+\langle p_1,\mu(s,x(s),u(s)\rangle+\langle p_2,\sigma_1(s,x(s),u(s),\sigma_2(s,x(s),u(s))\rangle.
	\end{align}
	Furthermore, assume $\tilde x(s)$ is the optimal probability to score a goal corresponding to the stubbornness $\tilde u(s)$. Then
	\begin{align}\label{2.47}
	dp(t)&=\mathcal{P}_x(s,\tilde x(s),\tilde u(s),p(s),q(s))\ ds+ q(s)\ dW_s,
	\end{align}
	has a mild unique solution.
	\label{th1}
\end{thm}

\begin{proof}
	Now we have a new BPPSDE. Lemmas 5.2 and 5.3 in \cite{du2013}, and proposition \ref{prop1}  directly imply this result.
\end{proof}

\subsection{Payoff function.}

In this section we formally define the expected payoff function of a single soccer player. Let $\E$ be the expectation with respect to the probability measure $\mathbb P$ and $\E_0\{.|\mathcal F_t\}=\E\{.|\mathcal F_t,x(0)=0\}$ is a conditional expectation based on fixed random initial goal dynamics $x(0)$, and the stubbornness $u(s)$ is adapted to filtration $\mathcal F_t$. Define an expected payoff function 
\begin{equation}\label{epayoff}
J(u):=\E_0\left\{\int_0^t\exp(-rs)\pi(s,x(s),u(s))ds+M\left(x(t)\right)\bigg|\mathcal F_0\right\},
\end{equation}
where $M\left(x(t)\right)$ is the terminal condition, and an individual soccer player's payoff function can be defined as 
\begin{equation}\label{payoff}
\pi(s,x(s),u(s)):=\left[\theta +\sum_{i=1}^3\a_i\right]x(s)-\frac{c(u(s))^2}{(r-\bar\mu)\sqrt{x(s)}},
\end{equation}
where $\theta>0$ is risk of injury, $\a_1$ is constant assist rate at time $s$, $\a_2$ is constant pass accuracy at time $s$, $\a_3$ is the dribbling skill of a player, $c>0$ represents a constant marginal cost of the player comes with sacrifices made by that player to be an eligible member of the team, and $\bar\mu$ is the average drift coefficient of the BPPSDE expressed in Equation \eqref{0}. Moreover $\frac{c(u(s))^2}{(r-\bar\mu)\sqrt{x(s)}}$ is the performance cost of a player, and $M(x(t))=\omega\exp(-rt)\sqrt{x(t)}$ is the terminal bonus offered to each player, where $\omega$ is some positive constant \citep{yeung2006cooperative}. The soccer player's objective is to find $J(u^*)=\sup_u J(u)$. The payoff function $\pi:I\times\mathcal X\times\mathcal U\mapsto \mathbb R$ for all $\mathcal X\subseteq\mathbb R$ and $\mathcal U\subseteq\mathbb R$, has a local maximum at $u^*\in\mathcal U$ if there exists a finite number $\epsilon>0$ so that for some  $u\in(u^*-\epsilon,u^*+\epsilon)\subset\mathcal U$, $J(u^*)\geq J(u)$.

\begin{lem}
	If the payoff function $J(.)$ is differentiable in the open functional space $\mathcal U$, and if it attains a local maximum at $u^*\in\mathcal U$ then $dJ/du=0$.
\end{lem}
\begin{proof}
	Let $J$ has a local maximum at $u^*\in\mathcal U$ if there exists a finite number $\epsilon>0$ so that for some  $u\in(u^*-\epsilon,u^*+\epsilon)\subset\mathcal U$, $J(u^*)\geq J(u)$. The first order total  derivative of $J$ w.r.t. $u$ is
	\[
	\frac{d}{du}J(u)=\lim_{u\ra u^*}\frac{J(u)-J(u^*)}{u-u^*}.
	\]
	Since $J(u^*)$ is maximum the numerator of the limit is never positive, but the denominator is positive in the cases where $u>u^*$ and negative for $u<u^*$. As we assume $J(.)$ is differentiable at $u^*$, the left and right limits exists and they are equal. This occurs if $dJ/du=0$. this completes the proof.
\end{proof}	

We assume $J(.)$ is continuously differentiable (i.e., smooth) w.r.t. $u$ in $(u^*-\epsilon,u^*+\epsilon)\subset\mathcal U,\ \ \forall\epsilon>0$. Define $u-u^*=\epsilon\kappa$. For $\epsilon\ra 0$ a second order Taylor series expansion implies
\[
J(u)=j(u^*)+\epsilon\kappa\frac{d}{du}J(u)+\frac{\epsilon^2}{2!}\kappa^2\frac{d^2}{du^2}J(u)+\mathcal O(\epsilon^3).
\]
If $dJ/du\neq0$ and $\epsilon\downarrow 0$, the sign of $J(u)-J(u^*)$ is unchanged in $(u^*-\epsilon,u^*+\epsilon)$, such that $\kappa dJ/du$ would have the same sign for every $\kappa$. It is trivial to understand that $\kappa$ can be non-zero. This implies $\kappa dJ/du$ can be non-zero too. Hence, $\kappa dJ/du=0$. Therefore, Taylor series expansion yields that the sign of the difference in stubbornness of a player is that of the quadratic term. If $d^2J/du^2<0$ then $J(u)$ attains a local maximum.

\section{Computation of the optimal stubbornness.}

In this section we are going construct a stochastic Lagrangian based on the system consists on Equations \eqref{epayoff} and \eqref{0}. In soccer, designing a player stubbornness that rewards the past performance is, in many ways, a much simpler task than designing one to predict future performance. This exercise boils down to assigning values to the stubbornness in decision making of  a player during a match or matches.  Therefore, th objective of a player is to maximize \eqref{epayoff} subject to the BPPSDE \eqref{0}. From Equation (45) of Ewald et al. (2024) \cite{ewald2024adaptation} for a player,  the stochastic Lagrangian at time $s\in[0,t]$ is defined as
\begin{multline}
\hat{\mathcal{L}}\left(s,{x},\lambda,u\right)=\E\biggr\{\int_0^t \bigg\{\exp(-rs)\pi(s,x(s),u(s))ds+M\left(x(t)\right)\biggr\}ds\\
+\int_0^t\left[x(s)-x_0-\int_0^s[\mu[\nu,x(\nu),u(\nu),\sigma_2[\nu,x(\nu),u(\nu)]]d\nu-\sigma[\nu,x(\nu),u(\nu)]dB(\nu)]\right] d\lambda(s)\biggr\},
\end{multline}

where $\lambda(s)$ is the Lagrangian multiplier.

\begin{prop}\label{p0}
	For a player if $X_0=\{ x(s), s\in[0,t]\}$ is a goal dynamics then, the stubbornness as a feedback Nash equilibrium $\big\{u^{*}(s,x)\in\mathcal U\big\}$ would be the solution of the following equation
	\begin{align}\label{11}
	\mbox{$\frac{\partial}{\partial u}$}f(s, x,u) \left[\mbox{$\frac{\partial^2}{\partial (x)^2}$}f(s, x,u)\right]^2=2\mbox{$\frac{\partial}{\partial x}$}f(s, x,u) \mbox{$\frac{\partial^2}{\partial x\partial u}$}f(s, x,u),
	\end{align}
	where for an It\^o process $h(s,x)\in \mathcal [0,t]\times\mathbb R$
	\begin{align}\label{12}
	f(s, x,u)&=\exp(-rs)\pi(s,x(s),u(s))ds+M\left(x(t)\right)+h(s,x)d\lambda(s)+\left[\mbox{$\frac{\partial h(s,x)}{\partial s}$}d\lambda(s)+\mbox{$\frac{d\lambda(s)}{d s}$}h(s,x)\right]\notag\\
	&\hspace{.25cm}+\mbox{$\frac{\partial h(s,x)}{\partial x}$}\mu\left[s,x,u,\sigma_2(s,x,u)\right]d\lambda(s)+\mbox{$\frac{1}{2}$}\left[\sigma\left[s,x,u\right]\right]^2\mbox{$\frac{\partial^2 h(s,x)}{\partial (x)^2}$}d\lambda(s).
	\end{align}
\end{prop}

\begin{proof}
	The Euclidean action function of a player can be represented as 
	\begin{align}
	\mathcal A_{0,t}(x)&=\int_0^t\E_s\bigg\{\exp(-rs)\pi(s,x(s),u(s))ds+M\left(x(t)\right)\notag\\
	&\hspace{1cm}+\bigg[x(s)-x_0-\mu\left[s,x,u,\sigma_2(s,x,u)\right]ds-\sigma\left[s,x,u\right]dB(\nu)\bigg]d\lambda(s)\bigg\},\notag
	\end{align}
	where $E_s$ is the conditional expectation on goal dynamics $x(s)$ at the beginning of time $s$. For all $\varepsilon>0$, and the normalizing constant $L_\varepsilon>0$ , define a transitional probability in small time interval as
	\begin{align}\label{w16}
	\Psi_{s,s+\varepsilon}(x)&:=\frac{1}{L_\varepsilon} \int_{\mathbb{R}} \exp\biggr\{-\varepsilon  \mathcal{A}_{s,s+\varepsilon}(x)\biggr\} \Psi_s(x) dx(s),
	\end{align}	
	for $\epsilon\downarrow 0$ and $\Psi_s(x)$ is the value of the transition probability at $s$ and goal dynamics $x(s)$ with the initial condition $\Psi_0(x)=\Psi_0$.
	
	For continuous time interval $[s,\tau]$ where $\tau=s+\varepsilon$  the stochastic Lagrangian is
	\begin{align}\label{action}
	\mathcal{A}_{s,\tau}(x)&= \int_{s}^{\tau}\ \E_s\ \biggr\{\exp(-r\nu)\pi(\nu,x(\nu),u(\nu)) d\nu\notag\\
	&\hspace{1cm}+\bigg[x(\nu)-x_0-\mu\left[\nu,x,u,\sigma_2(\nu,x,u)\right]d\nu-\sigma\left[\nu,x,u\right]dB(\nu)\bigg]d\lambda(\nu)\bigg\},
	\end{align}
	with the constant initial condition $x(0)=x_0$.	This conditional expectation is valid when the stubbornness $u(\nu)$ of a player's goal dynamics is determined at time $\nu$ such that all other players' goal dynamics are given. The evolution takes place as the action function is stationary. Therefore, the conditional expectation with respect to time only depends on the expectation of initial time point of interval $[s,\tau]$.
	
	Fubini's Theorem implies,
	\begin{align}\label{action5}
	\mathcal{A}_{s,\tau}(x)&= \E_s\ \bigg\{ \int_{s}^{\tau}\exp(-r\nu)\pi(\nu,x(\nu),u(\nu)) d\nu\notag\\
	&\hspace{1cm}+\bigg[x(\nu)-x_0-\mu\left[\nu,x,u,\sigma_2(\nu,x,u)\right]d\nu-\sigma\left[\nu,x,u\right]dB(\nu)\bigg]d\lambda(\nu) \bigg\}.
	\end{align}
	By It\^o's Theorem \citep{oksendal2013stochastic} there exists a function $h[\nu,x(\nu)]\in C^2([0,\infty)\times\mathbb{R})$ such that  $Y(\nu)=h[\nu,x(\nu)]$ where $Y(\nu)$ is an It\^o process.
	
	Assuming 
	\[
	h[\nu+\Delta \nu,x(\nu)+\Delta x(\nu)]= x(\nu)-x_0-\mu\left[\nu,x,u,\sigma_2(\nu,x,u)\right]d\nu-\sigma\left[\nu,x,u\right]dB(\nu),
	\]
	Equation (\ref{action5}) implies,
	\begin{align}\label{action6}
	\mathcal{A}_{s,\tau}(x)&=\E_s \bigg\{ \int_{s}^{\tau}\exp(-r\nu)\pi(\nu,x(\nu),u(\nu)) d\nu+ h\left[\nu+\Delta \nu,x(\nu)+\Delta x(\nu)\right]d\lambda(\nu)\bigg\}.
	\end{align} 
	
	It\^o's Lemma implies,
	\begin{align}\label{action7}
	\varepsilon\mathcal{A}_{s,\tau}(x)&= \E_s \bigg\{\varepsilon \exp(-rs)\pi(s,x(s),u(s))+ \varepsilon h[s,x(s)]d\lambda(s)+ \varepsilon h_s[s,x(s)]d\lambda(s) \notag\\
	&\hspace{.25cm}+\varepsilon h_x[s,x(s)]\mu\left[s,x(s),u(s),\sigma_2(s,x(s),u(s))\right]d\lambda(s) +\varepsilon h_x[s,x(s)]\sigma\left[s,x(s),u(s)\right]d\lambda(s) dB(s)\notag\\
	&\hspace{.5cm}+\mbox{$\frac{1}{2}$}\varepsilon\left(\sigma\left[s,x(s),u(s)\right]\right)^2h_{xx}[s,x(s)]d\lambda(s)+o(\varepsilon)\bigg\},
	\end{align}
	where $h_s=\frac{\partial}{\partial s} h$, $h_x=\frac{\partial}{\partial x} h$ and $h_{xx}=\frac{\partial^2}{\partial (x)^2} h$, and we use the condition $[d x(s)]^2\approx\varepsilon$ with $d x(s)\approx\varepsilon\mu[s,x(s),u(s),\sigma_2(s,x(s),u(s))]+\sigma[s,x(s),u(s)]dB(s)$.
	
	 We use It\^o Lemma and a similar approximation to approximate the integral. With $\varepsilon\downarrow 0$, dividing throughout by  $\varepsilon$ and taking the conditional expectation yields,
	\begin{align}\label{action8}
	\varepsilon\mathcal{A}_{s,\tau}(x)&= \E_s \bigg\{\varepsilon \exp(-rs)\pi(s,x(s),u(s))+\varepsilon h[s,x(s)]d\lambda(s)+ \varepsilon h_s[s,x(s)]d\lambda(s)\notag\\
	&\hspace{.25cm}+\varepsilon h_x[s,x(s)]\mu\left[s,x(s),u(s),\sigma_2(s,x(s),u(s))\right]d\lambda(s)\notag\\
	&\hspace{.5cm}+\mbox{$\frac{1}{2}$}\varepsilon\sigma^{2}\left[s,x(s),u(s)\right]h_{xx}[s,x(s)]d\lambda(s)+o(1)\bigg\},
	\end{align}
	since $\E_s[dB(s)]=0$ and $\E_s[o(\varepsilon)]/\varepsilon\ra 0$ for all $\varepsilon\downarrow 0$. For $\varepsilon\downarrow 0$ denote a transition probability at $s$ as $\Psi_s(x)$. Hence, using Equation (\ref{w16}), the transition function yields
	\begin{multline}\label{action9}
	\Psi_{s,\tau}(x)=\frac{1}{L_\epsilon^i}\int_{\mathbb{R}} \exp\biggr\{-\varepsilon \big[\exp(-rs)\pi(s,x(s),u(s))+h[s,x(s)]d\lambda(s)\\
	 +h_s[s,x(s)]d\lambda(s) +h_x[s,x(s)]\mu\left[s,x(s),u(s),\sigma_2(s,x(s),u(s))\right]d\lambda(s)\\
	+\mbox{$\frac{1}{2}$}\left(\sigma\left[s,x(s),u(s)\right]\right)^2h_{xx}[s,x(s)]d\lambda(s)\big]\biggr\} \Psi_s(x) dx(s)+o(\varepsilon^{1/2}).
	\end{multline}
	
	Since $\varepsilon\downarrow 0$, first order Taylor series expansion on the left hand side of Equation (\ref{action9}) yields
	\begin{multline}\label{action10}
	\Psi_{s}(x)+\varepsilon  \frac{\partial \Psi_{s}(x) }{\partial s}+o(\varepsilon)=\frac{1}{L_\varepsilon}\int_{\mathbb{R}} \exp\biggr\{-\varepsilon \big[\exp(-rs)\pi(s,x(s),u(s))+h[s,x(s)]d\lambda(s) \\
	+h_s[s,x(s)]d\lambda(s)+h_x[s,x(s)]\mu\left[s,x(s),u(s),\sigma_2(s,x(s),u(s))\right]d\lambda(s)\\
	+\mbox{$\frac{1}{2}$}\left(\sigma\left[s,x(s),u(s)\right]\right)^2h_{xx}[s,x(s)]d\lambda(s)\big]\biggr\} \Psi_s(x) dx(s)+o(\varepsilon^{1/2}).
	\end{multline}
	For any given $s$ and $\tau$ define $x(s)-x(\tau):=\xi$ such that $x(s)=x(\tau)+\xi$. For the instance where $\xi^i$ is not around zero, for a positive number $\eta<\infty$ assume $|\xi|\leq\sqrt{\frac{\eta\varepsilon}{x(s)}}$ such that for $\varepsilon\downarrow 0$, $\xi$ attains smaller values and the goal dynamics $0<x(s)\leq\eta\varepsilon/(\xi)^2$. Thus,
	\begin{multline*}
	\Psi_{s}(x)+\varepsilon\frac{\partial \Psi_{s}(x)}{\partial s}=\frac{1}{L_\epsilon}\int_{\mathbb{R}} \left[\Psi_{s}(x)+\xi\frac{\partial \Psi_{is}(x)}{\partial x}+o(\epsilon)\right]\\
	\times \exp\biggr\{-\varepsilon \big[\exp(-rs)\pi(s,x(s),u(s))+h[s,x(s)]d\lambda(s)+h_x[s,x(s)]\mu\left[s,x(s),u(s),\sigma_2(s,x(s),u(s))\right]d\lambda(s)\\
	+\mbox{$\frac{1}{2}$}\left(\sigma\left[s,x(s),u(s)\right]\right)^2h_{xx}[s,x(s)]d\lambda(s)\big]\biggr\} d\xi+o(\varepsilon^{1/2}).
	\end{multline*}
	
	Before solving for Gaussian integral of the each term of the right hand side of the above Equation define a $C^2$ function 
	\begin{align*}
	f[s,\xi,\lambda(s),u(s)]&=\exp(-rs)\pi(s,x(s)+\xi,u(s))+h[s,x(s)+\xi]d\lambda(s) +h_s[s,x(s)+\xi]d\lambda(s)\notag\\
	&\hspace{.25cm}+h_x[s,x(s)+\xi]\mu\left[s,x(s)+\xi,u(s),\sigma_2(s,x(s)+\xi,u(s))\right]d\lambda(s)\\
	&\hspace{1cm}+\mbox{$\frac{1}{2}$}\sigma^{2}\left[s,x(s)+\xi,u(s)\right]h_{xx}[s,x(s)+\xi]d\lambda(s)+o(1).
	\end{align*}
	 Hence,
	\begin{align}\label{action13}
	\Psi_{s}(x)+\varepsilon \frac{\partial \Psi_{s}(x) }{\partial s}&=\Psi_{s}(x)\frac{1}{L_\epsilon}\int_{\mathbb{R}}\exp\left\{-\varepsilon f[s,\xi,\lambda(s),u(s)] \right\}d\xi\notag\\
	&+\frac{\partial \Psi_{s}(x)}{\partial x}\frac{1}{L_\epsilon}\int_{\mathbb{R}}\xi\exp\left\{-\varepsilon f[s,\xi,\lambda(s),u(s)] \right\}d\xi+o(\varepsilon^{1/2}).
	\end{align}
	
	After taking $\varepsilon\downarrow 0$, $\Delta u\downarrow0$ and a Taylor series expansion with respect to $x$ of $f[s,\xi,\lambda(s),u(s)]$ yields, 
	\begin{align*}
	f[s,\xi,\lambda(s),u(s)]&=f[s, x(\tau),\lambda(s),u(s)]+f_x[s, x(\tau),\lambda(s),\gamma,u(s)][\xi-x(\tau)]\notag\\
	&\hspace{1cm}+\mbox{$\frac{1}{2}$}f_{xx}[s, x(\tau),\lambda(s),u(s)][\xi-x(\tau)]^2+o(\varepsilon).
	\end{align*}
	Define $y:=\xi-x(\tau)$ so that $ d\xi=dy$. First integral on the right hand side of Equation (\ref{action13}) yields
	\begin{align}\label{action14}
	&\int_{\mathbb{R}} \exp\big\{-\varepsilon f[s,\xi,\lambda(s),u(s)]\} d\xi\notag\\
	&=\exp\big\{-\varepsilon f[s, x(\tau),\lambda(s),u(s)]\big\}\notag\\
	&\hspace{1cm}\int_{\mathbb{R}} \exp\biggr\{-\varepsilon \biggr[f_x[s, x(\tau),\lambda(s),u(s)]y+\mbox{$\frac{1}{2}$}f_{xx}[s,x(\tau),\lambda(s),u(s)](y^i)^2\biggr]\biggr\} dy.
	\end{align}
	Assuming  $a=\frac{1}{2} f_{xx}[s, x(\tau),\lambda(s),u(s)]$ and $b=f_x[s, x(\tau),\lambda(s),u(s)]$ the argument of the exponential function in Equation (\ref{action14}) becomes,
	\begin{align}\label{action15}
	a(y)^2+by&=a\left[(y)^2+\frac{b}{a}y\right]=a\left(y+\frac{b}{2a}y\right)^2-\frac{(b)^2}{4(a)^2}.
	\end{align}
	
	Therefore,
	\begin{align}\label{action16}
	&\exp\bigg\{-\varepsilon f[s, x(\tau),\lambda(s),u(s)]\bigg\}\int_{\mathbb{R}} \exp\big\{-\varepsilon [a(y)^2+by]\big\}dy\notag\\
	&=\exp\left\{\varepsilon \left[\frac{b^2}{4a^2}-f[s, x(\tau),\lambda(s),u(s)]\right]\right\}\int_{\mathbb{R}} \exp\left\{-\left[\varepsilon a\left(y+\frac{b}{2a}y\right)^2\right]\right\} dy\notag\\
	&=\sqrt{\frac{\pi}{\varepsilon a}}\exp\left\{\varepsilon \left[\frac{b^2}{4a^2}-f[s, x(\tau),\lambda(s),u(s)]\right]\right\},
	\end{align}
	and
	\begin{align}\label{action17}
	&\Psi_{s}(x) \frac{1}{L_\varepsilon} \int_{\mathbb{R}} \exp\big\{-\varepsilon f[s,\xi,\lambda(s),u(s)]\} d\xi\notag\\
	 &=\Psi_{s}(x) \frac{1}{L_\varepsilon} \sqrt{\frac{\pi}{\varepsilon a}}\exp\left\{\varepsilon \left[\frac{b^2}{4a^2}-f[s, x(\tau),\lambda(s),u(s)]\right]\right\}. 
	\end{align}
	Substituting $\xi=x(\tau)+y$ second integrand of the right hand side of Equation (\ref{action13}) yields
	\begin{align}\label{action18}
	& \int_{\mathbb{R}} \xi \exp\left[-\varepsilon \{f[s,\xi,\lambda(s),u(s)]\}\right] d\xi\notag\\
	&=\exp\{-\varepsilon f[s, x(\tau),\lambda(s),u(s)]\}\int_{\mathbb{R}} [x(\tau)+y] \exp\left[-\varepsilon \left[ay^2+by\right]\right] dy\notag\\
	&=\exp\left\{\varepsilon \left[\frac{b^2}{4a^2}-f[s, x(\tau),\lambda(s),u(s)]\right]\right\} \biggr[x(\tau)\sqrt{\frac{\pi}{\varepsilon a}}\notag\\
	&\hspace{1cm}+\int_{\mathbb{R}} y \exp\left\{-\varepsilon \left[a\left(y+\frac{b}{2a}y\right)^2\right]\right\} dy\biggr].
	\end{align}
	
	Substituting $k=y+b/(2a)$ in Equation (\ref{action18}) yields,
	\begin{align}\label{action19}
	&\exp\left\{\varepsilon \left[\frac{b^2}{4a^2}-f^i[s, x(\tau),\lambda(s),u(s)]\right]\right\} \biggr[x(\tau)\sqrt{\frac{\pi}{\varepsilon a}}+\int_{\mathbb{R}} \left(k-\frac{b}{2a}\right) \exp[-a\varepsilon k^2] dk\biggr]\notag\\
	&=\exp\left\{\varepsilon \left[\frac{b^2}{4a^2}-f[s, x(\tau),\lambda(s),u(s)]\right]\right\} \biggr[x(\tau)-\frac{b}{2a}\biggr]\sqrt{\frac{\pi}{\varepsilon a}}.
	\end{align}
	Hence,
	\begin{align}\label{action20}
	&\frac{1}{L_\varepsilon}\frac{\partial \Psi_{s}(x)}{\partial x}\int_{\mathbb{R}} \xi \exp\left[-\varepsilon f[s,\xi,\lambda(s),u(s)]\right] d\xi\notag\\
	&=\frac{1}{L_\varepsilon}\frac{\partial \Psi_{s}(x)}{\partial x} \exp\left\{\varepsilon \left[\frac{b^2}{4a^2}-f[s, x(\tau),\lambda(s),u(s)]\right]\right\} \biggr[x(\tau)-\frac{b}{2a}\biggr]\sqrt{\frac{\pi}{\varepsilon a}}.
	\end{align}
	Plugging in Equations (\ref{action17}), and (\ref{action20})  into Equation (\ref{action13}) implies,
	\begin{align}\label{action24}
	&\Psi_{s}(x)+\varepsilon \frac{\partial \Psi_{s}(x)}{\partial s}\notag\\
	&=\frac{1}{L_\varepsilon} \sqrt{\frac{\pi}{\varepsilon a}}\Psi_{s}(x) \exp\left\{\varepsilon \left[\frac{b^2}{4a^2}-f[s, x(\tau),\lambda(s),u(s)]\right]\right\}\notag\\
	&+\frac{1}{L_\varepsilon}\frac{\partial \Psi_{s}(x)}{\partial x} \sqrt{\frac{\pi}{\varepsilon a}} \exp\left\{\varepsilon \left[\frac{b^2}{4a^2}-f[s, x(\tau),\lambda(s),u(s)]\right]\right\} \biggr[x(\tau)-\frac{b}{2a}\biggr]+o(\varepsilon^{1/2}).
	\end{align}
	Let $f^i$ be in Schwartz space. This leads to derivatives are rapidly falling  and  furthermore, assuming $0<|b|\leq\eta\varepsilon$, $0<|a|\leq\mbox{$\frac{1}{2}$}[1-(\xi)^{-2}]^{-1}$ and $x(s)-x(\tau)=\xi$ yields,
	\begin{align*}
	x(\tau)-\frac{b}{2a}=x(s)-\xi-\frac{b}{2a}=x(s)-\frac{b}{2a},\ \forall\ \xi\downarrow 0,
	\end{align*}
	such that 
	\begin{align*}
	\bigg|x(s)-\frac{b}{2a}\bigg|=\bigg|\frac{\eta\varepsilon}{(\xi)^2}-\eta\varepsilon\left[1-\frac{1}{(\xi)^2}\right]\bigg|\leq\eta\varepsilon.
	\end{align*}
	
	Therefore, the Wick rotated Schr\"odinger-type Equation for the player is,
	\begin{align}\label{action25.4}
	\frac{\partial \Psi_{s}(x)}{\partial s}&=\left[\frac{b^2}{4a^2}-f[s, x(\tau),\lambda(s),u(s)]\right]\Psi_{s}(x).
	\end{align}
	Differentiating the Equation (\ref{action25.4}) with respect to stubbornness yields 
	\begin{align}\label{w18}
	\left\{\frac{2f_x}{f_{xx}}\left[\frac{f_{xx}f_{xu}-f_xf_{xxu}}{(f_{xx})^2}\right]-f_u\right\}\Psi_{s}(x)=0,
	\end{align}
	where $f_x=\frac{\partial}{\partial x} f$, $f_{xx}=\frac{\partial^2}{\partial (x)^2} f$, $f_{xu}=\frac{\partial^2}{\partial x\partial u} f$ and $f_{xxu}=\frac{\partial^3}{\partial (x)^2\partial u} f=0$. Thus, optimal feedback stubbornness of a player in stochastic goal dynamics is represented as $u^{*}(s,x)$ and is found by setting Equation (\ref{w18}) equal to zero. Hence, $u^{*}(s,x)$ is the solution of the following Equation
	\begin{align}\label{w21}
	f_u (f_{xx})^2=2f_x f_{xu}.
	\end{align}
	This completes the proof.
\end{proof}

\begin{remark}
	The central idea of Proposition \ref{p0} is to choose $h$ appropriately. Therefore, one natural candidate should be a function of the integrating factor of the stochastic goal dynamics represented in Equation \eqref{0}.
\end{remark}

To demonstrate the preceding proposition, we present a detailed example to identify an optimal stubbornness under this environment. Consider a player has to maximize the expected payoff of expressed in Equation \eqref{payoff}
\begin{equation}\label{e0}
J(u):=\E_0\left\{\int_0^t\exp(-rs)\left\{\left[\theta +\sum_{i=1}^3\a_i\right]x(s)-\frac{c(u(s))^2}{(r-\bar\mu)\sqrt{x(s)}}\right\}ds+\omega\exp(-rt)\sqrt{x(t)}\bigg|\mathcal F_0\right\},
\end{equation}
subject to the goal dynamics represented by a BPPSDE
\begin{equation}\label{e1}
dx(s)=\left[a\sqrt{x(s)}-\sigma_2x(s)-u(s)\right]ds+\left[\sigma_1-\sigma_2x(s)\right]dB(s),
\end{equation}
where $a$ is a constant, $\sigma_1$ and $\sigma_2$ are constant volatilities, and the diffusion component is $\left[\sigma_1-\sigma_2x(s)\right]$. We are going to implement Proposition \ref{p0} to determine the optimal stubbornness. By this problem
\begin{multline}\label{e2}
f(s,x,u)=\exp(-rs)\left\{\left[\theta +\sum_{i=1}^3\a_i\right]x(s)-\frac{c(u(s))^2}{(r-\bar\mu)\sqrt{x(s)}}\right\}+\omega\exp(-rt)\sqrt{x(t)}\\
+h(s,x)d\lambda(s)+\left[\mbox{$\frac{\partial h(s,x)}{\partial s}$}d\lambda(s)+\mbox{$\frac{d\lambda(s)}{d s}$}h(s,x)\right]+\mbox{$\frac{\partial h(s,x)}{\partial x}$}\left[a\sqrt{x(s)}-\sigma_2x(s)-u(s)\right]d\lambda(s)\\
+\frac{1}{2}\left[\sigma_1-\sigma_2x(s)\right]^2\mbox{$\frac{\partial^2 h(s,x)}{\partial (x)^2}$}d\lambda(s).
\end{multline}
In Equation \eqref{e2} we treat the terminal condition as a constant hence, define $\omega\exp(-rt)\sqrt{x(t)}=\bar M$. Since the diffusion part $\left[\sigma_1-\sigma_2x(s)\right]$ of the SDE \eqref{e1} suggests that we can simplify the equation by focusing on this part. One common approach is to consider an exponential integrating factor to counterbalance the $-\sigma_2x(s)$ term in both the drift and diffusion terms. Therefore, the integrating factor $h(s,x)=\exp(\sigma_2 x(s))$. Therefore, Equation \eqref{e2} yields
\begin{multline}\label{e3}
f(s,x,u)=\exp(-rs)\left\{\left[\theta +\sum_{i=1}^3\a_i\right]x(s)-\frac{c(u(s))^2}{(r-\bar\mu)\sqrt{x(s)}}\right\}+\bar M\\
+\exp(\sigma_2x(s))d\lambda(s)+\mbox{$\frac{d\lambda(s)}{d s}$}\exp(\sigma_2x(s))+\sigma_2\exp(\sigma_2x(s))\left[a\sqrt{x(s)}-\sigma_2x(s)-u(s)\right]d\lambda(s)\\
+\frac{1}{2}\left[\sigma_1-\sigma_2x(s)\right]^2(\sigma_2)^2\exp(\sigma_2x(s)).
\end{multline}
Hence,
\begin{equation}\label{e4}
\frac{\partial }{\partial u}f(s,x,u)=-\exp(-rs)\frac{2cu(s)}{(r-\bar\mu)\sqrt{x(s)}}-\sigma_2\exp(\sigma_2x(s))d\lambda(s),
\end{equation}
\begin{multline}\label{e5}
\frac{\partial f(s, x, u)}{\partial x(s)} = \exp(-rs)\left\{ \left[ \theta + \sum_{i=1}^3 \alpha_i \right] -   \frac{c(u(s))^2}{2(r - \bar{\mu}) x(s)^{3/2}}\right\}\\
+ \sigma_2 \exp(\sigma_2 x(s)) \bigg[d\lambda(s) + \frac{d\lambda(s)}{ds}
+  \left( \frac{a}{2 \sqrt{x(s)}} - \sigma_2 \right) d\lambda(s)\\
+ \sigma_2  \left[ a \sqrt{x(s)} - \sigma_2 x(s) - u(s) \right] d\lambda(s)\bigg]
- \left[ \sigma_1 - \sigma_2 x(s) \right] (\sigma_2)^3 \exp(\sigma_2 x(s))\\
+ \frac{1}{2} \left[\sigma_1 - \sigma_2 x(s) \right]^2 (\sigma_2)^3 \exp(\sigma_2 x(s)),
\end{multline}
\begin{multline}\label{e6}
\frac{\partial^2 f(s, x, u)}{\partial x(s)^2} = \exp(-rs) \frac{15c(u(s))^2}{4(r - \bar{\mu}) x(s)^{5/2}}+ (\sigma_2)^2 \exp(\sigma_2 x(s)) d\lambda(s)+ (\sigma_2)^2 \frac{d\lambda(s)}{ds} \exp(\sigma_2 x(s))\\
+ (\sigma_2)^2 \exp(\sigma_2 x(s)) \left( \frac{a}{2 \sqrt{x(s)}} - \sigma_2 \right) d\lambda(s)
- \sigma_2 \exp(\sigma_2 x(s)) \frac{3a}{4 x(s)^{5/2}} d\lambda(s)\\
+ (\sigma_2)^3 \exp(\sigma_2 x(s)) \left[ a \sqrt{x(s)} - \sigma_2 x(s) - u(s) \right] d\lambda(s)\\
+ (\sigma_2)^2 \exp(\sigma_2 x(s)) \frac{a}{2 \sqrt{x(s)}} d\lambda(s)
- (\sigma_2)^2 \exp(\sigma_2 x(s)) \frac{a}{4 x(s)^{3/2}} d\lambda(s)\\
+ (\sigma_2)^4 \exp(\sigma_2 x(s))
- \left[ \sigma_1 - \sigma_2 x(s) \right] (\sigma_2)^5 \exp(\sigma_2 x(s))\\
+ \frac{1}{2} \left[ \sigma_1 - \sigma_2 x(s) \right]^2 (\sigma_2)^4 \exp(\sigma_2 x(s)),
\end{multline}
and
\begin{equation}\label{e7}
\frac{\partial^2 f(s, x, u)}{\partial x \partial u} = -\frac{c \exp(-rs) u(s)}{(r - \bar{\mu}) x(s)^{3/2}}.
\end{equation}
Equation \eqref{w18} yields,
\begin{multline}\label{e8}
\left( -\exp(-rs) \frac{2 c u(s)}{(r - \bar{\mu}) \sqrt{x(s)}} - A_1 \right) \left( \exp(-rs) \frac{15 c (u(s))^2}{4 (r - \bar{\mu}) x(s)^{5/2}} - (\sigma_2)^3 \exp(\sigma_2 x(s)) u(s) d\lambda(s) + A_3 \right)^2\\
= 2 \left( -\exp(-rs) \frac{c (u(s))^2}{2 (r - \bar{\mu}) x(s)^{3/2}} - (\sigma_2)^2 \exp(\sigma_2 x(s)) u(s) d\lambda(s) + A_2 \right) \left( -\frac{c \exp(-rs) u(s)}{(r - \bar{\mu}) x(s)^{3/2}} \right),
\end{multline}
where
\begin{align*}
A_1&=\sigma_2\exp(\sigma_2x(s))d\lambda(s),\\
A_2&=\exp(-rs)\left[\theta+\sum_{i=1}^3\a_i\right]+\sigma_2 \exp(\sigma_2 x(s)) \bigg[d\lambda(s) + \frac{d\lambda(s)}{ds}
+  \left( \frac{a}{2 \sqrt{x(s)}} - \sigma_2 \right) d\lambda(s)\\
&\hspace{1cm}+ \sigma_2  \left[ a \sqrt{x(s)} - \sigma_2 x(s) \right] d\lambda(s)\bigg]
- \left[ \sigma_1 - \sigma_2 x(s) \right] (\sigma_2)^3 \exp(\sigma_2 x(s))\\
&\hspace{2cm}+ \frac{1}{2} \left[\sigma_1 - \sigma_2 x(s) \right]^2 (\sigma_2)^3 \exp(\sigma_2 x(s)),\\
A_3&=(\sigma_2)^2 \exp(\sigma_2 x(s)) d\lambda(s)+ (\sigma_2)^2 \frac{d\lambda(s)}{ds} \exp(\sigma_2 x(s))\\
&\hspace{.5cm}+ (\sigma_2)^2 \exp(\sigma_2 x(s)) \left( \frac{a}{2 \sqrt{x(s)}} - \sigma_2 \right) d\lambda(s)
- \sigma_2 \exp(\sigma_2 x(s)) \frac{3a}{4 x(s)^{5/2}} d\lambda(s)\\
&\hspace{1cm}+ (\sigma_2)^3 \exp(\sigma_2 x(s)) \left[ a \sqrt{x(s)} - \sigma_2 x(s) \right] d\lambda(s)\\
&\hspace{1.5cm}+ (\sigma_2)^2 \exp(\sigma_2 x(s)) \frac{a}{2 \sqrt{x(s)}} d\lambda(s)
- (\sigma_2)^2 \exp(\sigma_2 x(s)) \frac{a}{4 x(s)^{3/2}} d\lambda(s)\\
&\hspace{2cm}+ (\sigma_2)^4 \exp(\sigma_2 x(s))
- \left[ \sigma_1 - \sigma_2 x(s) \right] (\sigma_2)^5 \exp(\sigma_2 x(s))\\
&\hspace{2.5cm}+ \frac{1}{2} \left[ \sigma_1 - \sigma_2 x(s) \right]^2 (\sigma_2)^4 \exp(\sigma_2 x(s)).
\end{align*}
Given the complexity of the terms (including exponents and products), the equation may result in a high-degree polynomial in 
$u(s)$, and explicit solution of stubbornness requires further simplification or assumptions to reduce the degree of the polynomial. Assume the effect of $d\lambda(s)$ is very small (i.e., $d\lambda(s)\ra 0$). This would remove $A_1$, $(\sigma_2)^3 \exp(\sigma_2 x(s)) u(s) d\lambda(s)$, and $(\sigma_2)^2 \exp(\sigma_2 x(s)) u(s) d\lambda(s)$ from Equation \eqref{e8}. Equation \eqref{e8} imples
\begin{multline}\label{e9}
\left( -\exp(-rs) \frac{2 c u(s)}{(r - \bar{\mu}) \sqrt{x(s)}} \right) \left( \exp(-rs) \frac{15 c (u(s))^2}{4 (r - \bar{\mu}) x(s)^{5/2}} + A_3 \right)^2\\
=2 \left( -\exp(-rs) \frac{c (u(s))^2}{2 (r - \bar{\mu}) x(s)^{3/2}} + A_2 \right) \left( -\frac{c \exp(-rs) u(s)}{(r - \bar{\mu}) x(s)^{3/2}} \right).
\end{multline}
Simplifying the left and the right hand sides of Equation \eqref{e9} we get
\[
-\exp(-2rs) \frac{2 c u(s)}{(r - \bar{\mu}) \sqrt{x(s)}} \left( \frac{15 c (u(s))^2}{4 (r - \bar{\mu}) x(s)^{5/2}} + A_3 \right)^2,
\]
and 
\[
\exp(-2rs) \frac{c^2 u(s)^3}{(r - \bar{\mu})^2 x(s)^{3}} - 2 A_2 \frac{c \exp(-rs) u(s)}{(r - \bar{\mu}) x(s)^{3/2}},
\]
respectively. Therefore, equating both sides and dividing by $\exp(-2rs)$ yield,
\begin{align}
&- \frac{2 c u(s)}{(r - \bar{\mu}) \sqrt{x(s)}} \left( \frac{15 c (u(s))^2}{4 (r - \bar{\mu}) x(s)^{5/2}} + A_3 \right)^2 =  \frac{c^2 u(s)^3}{(r - \bar{\mu})^2 x(s)^{3}} - 2 A_2 \frac{c  u(s)\exp(-3rs)}{(r - \bar{\mu}) x(s)^{3/2}}\notag\\
&-\frac{2 c u(s)}{(r - \bar{\mu}) \sqrt{x(s)}} \left( \frac{15 c (u(s))^2}{4 (r - \bar{\mu}) x(s)^{5/2}} + A_3 \right)^2 - \frac{c^2 u(s)^3}{(r - \bar{\mu})^2 x(s)^{3}} + 2 A_2 \frac{c u(s)\exp(-3rs)}{(r - \bar{\mu}) x(s)^{3/2}} = 0\notag\\
&u(s) \left( -\frac{2 c}{(r - \bar{\mu}) \sqrt{x(s)}} \left( \frac{15 c (u(s))^2}{4 (r - \bar{\mu}) x(s)^{5/2}} + A_3 \right)^2 - \frac{c^2 u(s)^2}{(r - \bar{\mu})^2 x(s)^{3}} + 2 A_2 \frac{c\exp(-3rs)}{(r - \bar{\mu}) x(s)^{3/2}} \right) = 0.
\end{align}
This gives one solution $u(s) = 0$ and we need to find the non-trivial solution by solving:
\[
-\frac{2 c}{(r - \bar{\mu}) \sqrt{x(s)}} \left( \frac{15 c (u(s))^2}{4 (r - \bar{\mu}) x(s)^{5/2}} + A_3 \right)^2 - \frac{c^2 u(s)^2}{(r - \bar{\mu})^2 x(s)^{3}} + 2 A_2 \frac{c\exp(-3rs)}{(r - \bar{\mu}) x(s)^{3/2}}=0.
\]
The above equation is quadratic in \( u(s) \). Let \( k_1 = -\frac{2 c}{(r - \bar{\mu}) \sqrt{x(s)}} \), \( k_2 = \frac{15 c}{4 (r - \bar{\mu}) x(s)^{5/2}}\), \( k_3 = \frac{c^2}{(r - \bar{\mu})^2 x(s)^{3}} \), and \( k_4 = 2 A_2 \frac{c\exp(-3rs)}{(r - \bar{\mu}) x(s)^{3/2}} \). We express the equation as $k_1k_2 u(s)^4 + (2k_1k_2A_3-k_3) u(s)^2 + (k_1A_3^2+k_4) = 0$. Define \(z(s):=u(s)^2\). Quadratic formula implies
\[
z^*(s)=\frac{1}{2}\left[\frac{1}{k_1k_2}\left(k_3-2k_1k_2\right)\pm\left[\frac{1}{(k_1k_2)^2}\left(k_3-2k_1k_2\right)^2-\frac{4}{k_1k_2}\left(k_1A_3^2-k_4\right)\right]^{1/2}\right].
\]
Therefore optimal stubbornness is,
\[
u^*(s)=\pm\left\{\frac{1}{2}\left[\frac{1}{k_1k_2}\left(k_3-2k_1k_2\right)\pm\left[\frac{1}{(k_1k_2)^2}\left(k_3-2k_1k_2\right)^2-\frac{4}{k_1k_2}\left(k_1A_3^2-k_4\right)\right]^{1/2}\right]\right\}^{1/2}.
\]
Since we assume the stubbornness is a non-negative function, we ignore the negative squared-root part.

\section{Conclusion.}
Over the past decade, research on modeling scores in soccer games has increasingly focused on dynamics to explain changes in team strengths over time. A crucial aspect of this involves evaluating the performance of all team players \citep{pramanik2024stc}. Consequently, a game-theoretic approach to determine optimal stubbornness has become essential. To compute this optimal stubbornness, we begin by constructing a stochastic Lagrangian based on the payoff function and the BPPSDE. We then apply a Euclidean path integral approach, derived from the Feynman action function \citep{feynman1948space}, over small continuous time intervals \citep{pramanik2025factors}. Through Taylor series expansion and Gaussian integral solutions, we derive a Wick-rotated Schr\"odinger-type equation. The analytical solution for optimal stubbornness is obtained by taking the first derivative with respect to stubbornness. This method simplifies challenges associated with the value function in the Hamiltonian-Jacobi-Bellman (HJB) equation. Moreover, under the BPPSDE framework, the path integral control approach performs more effectively than the HJB equation \citep{pramanik2025strategies}.

Considering the significant impact of stubbornness research on predicting match outcomes—not only for the betting industry but also for soccer clubs and analytics teams—there is strong motivation for researchers to explore this field. We believe many game theorists who are also soccer enthusiasts would welcome seeing this \emph{Moneyball effect} extend into soccer. Ultimately, these concepts could be applied beyond soccer to any sport where score predictions rely on strength dynamics \citep{pramanik2024motivation}.

\section*{Conflict of interest.}
The author declares that he has no conflict of interest.

\section*{Data Availability.}
The author declares that no data have been used in this paper.

\section*{Funding.}
This research received no external funding.

\section*{Appendix.}
\begin{prop}
	A stochastic differential equation can be represented as a Feynman path integration.
	\label{prop11}
\end{prop}
\begin{proof}
	For the simplistic case consider the Feynman action function has the form $$S(x)=-\left[\int_0^t\ \left\{ \pi(s,x(s),u(s))+\lambda\ g(s,x(s),\dot x(s))\right\}\ ds\right].$$
	 All the symbols have the same meaning as described in the main text. In this context we assume that the penalization function $\lambda g(s,x(s),\dot x(s))$ is a proxy of an SDE, and the payoff function is stable and it can be added to the drift part of the stochastic differential equation. Let the penalization function be of the form $$\frac{d x}{d s}=A(s,x(s),u(s))+\sigma(s,x(s),u(s))\ B_s.$$ Including a payoff function yields 
	\begin{align}\label{a0}
	\frac{d x}{d s}&=\pi(s,x(s),u(s))+A(s,x(s),u(s))+\sigma(s,x(s),u(s))\ B_s\notag\\i.e.\ \frac{d x}{d s}&=\mu(s,x(s),u(s))+ \sigma(s,x(s),u(s))\ B_s
	\end{align}
	where $\mu(s,x(s),u(s))=\pi(s,x(s),u(s))+A(s,x(s),u(s))$. Therefore, after using condition (\ref{a0}) our new general \emph{Langevin} form becomes,
	\begin{align}\label{a1}
	\frac{d x}{d s}&=\mu(s,x(s),u(s)) +\sigma(s,x(s),u(s)) B_s,
	\end{align}
	where $B_s$ is the \emph{white noise} with zero mean and variance $|s_i-s_j|$ for all ${i,j}\in\{0,1,...,n\}^2$ and $i\neq j$. Equation (\ref{a1}) can be written as,
	\begin{align}\label{a2}
	dx &= \mu(s,x(s),u(s))\ ds +\sigma(s,x(s),u(s))\ dW_s.
	\end{align}
	We further assume that, the functions $\mu(s,x(s),u(s))$ and $\sigma(s,x(s),u(s))$ obey all the properties of Ito stochastic differential equation. We want to derive a probability density function (PDF) for a soccer player's goal dynamics at time $s$ [i.e. $x(s)$]. The decentralized form of equation (\ref{a2}) with the small Ito interpretation of small time step $h$ is
	\begin{align}\label{a3}
	x_{j+1}-x_j&=\mu_j\ h+\sigma_j\ \omega_j\ \sqrt h,\ \forall j\in\{0,1,...,n\},
	\end{align}
	where initial time is zero, $x_j=x(0+jh)$, $T=nh$, $\mu_j=\mu(0+jh,x_j)$, $\sigma_j=\sigma(0+jh,x_j)$, $\omega_j$ is a normally distributed discrete random variable with $\langle w_j\rangle=0$ and $\langle w_iw_j\rangle=\Delta_{j,k}$. Chow (2015) \cite{chow2015} define $\Delta_{j,k}$ as the \emph{Kronecker delta} function.
	
	The conditional PDF of firm's output share can be written as,
	\begin{align}\label{a4}
	\Psi(x|\omega)&=\prod_{j=0}^n\ \delta(x_{j+1}-x_j-\mu_j h-\sigma_j\omega_j\sqrt h)
	\end{align}
	Probability in equation (\ref{a4}) is nothing but the delta \emph{Dirac} function constrained on the stochastic differential equation. We know the \emph{Fourier} transformation of delta Dirac function is,
	\begin{align}\label{a5}
	\delta(b_j)=\frac{1}{2\pi}\int_{\mathbb R}\ \exp\{-i\mathcal L_jb_j\}\ db_j
	\end{align}
	where $-i$ is the complex number. After putting the stochastic differential equation (\ref{a3})  into $\mathcal L_j$ we get,
	\begin{align}\label{a6}
	\Psi(x|\omega)&=\int_{\mathbb R}\ \prod_{j=0}^n\ \frac{1}{2\pi}\ e^{-i\sum_j\ b_j(x_{j+1}-x_j-\mu_j h-\sigma_j\omega_j\sqrt h)}\ db_j
	\end{align}
	\cite{chow2015} suggests that zero mean unit-variance \emph{Gaussian white noise} density function can be written as
	\begin{align}\label{a7}
	\Psi(\omega_j)=\frac{1}{\sqrt{2\pi}}\ e^{-(1/2)\omega^2_j}
	\end{align} 
	Combining equations (\ref{a6}) and (\ref{a7}) we get,
	\begin{align}\label{a8}
	\Psi(0,t,x_0,x_t)&=\int_{\mathbb R}\ \Psi(x|\omega)\ \prod_{j=0}^n\ \Psi(\omega_j)\ d\omega_j\notag\\ &=\int_{\mathbb R}\prod_{j=0}^n\frac{1}{2\pi}\ e^{-i\ \sum_j\ b_j(x_{j+1}-x_j-\mu_j h)}\ db_j\ \times \ \int_{\mathbb R}\prod_{j=0}^n\frac{1}{\sqrt{2\pi}}\ e^{-ib_j\sigma_j\omega_j\sqrt h}\ e^{-(1/2)\omega^2_j}\ d\omega_j\notag\\&=\int_{\mathbb R}\prod_{j=0}^n\frac{1}{2\pi}\ e^{-\sum_j\ (ib_j)\left(\left(\frac{x_{j+1}-x_j}{h}-\mu_j\right)h\right)+\sum_j\ (1/2)^2\ \sigma^2_j(ib_j)^2h}\ db_j
	\end{align}
	Since $h\ra0$, $n\ra\infty$ such that $T=nh$, Equation (\ref{a8}) yields 
	\begin{align}\label{a9}
	\Psi(0,t,x_0,x_t)&=\int_{\mathbb R}\ e^{-\int_0^t\ [\tilde x(s)(\dot x(s)-\mu(s,x(s),u(s)))-\frac{1}{2}\tilde x(s)^2 \sigma(s,x(s),u(s))]\ ds}\ \mathcal{D}_x
	\end{align}
	with a newly defined complex variable $ib_j\ra\tilde x$. Chow et al. (2015) \cite{chow2015} imply that although they use continuum notation for covariance, $x(s)$ does not need be differentiable and they interpret the action function by discreet definition. In the power of exponential integrand in Equation (\ref{a9}) there are two parts, one is real and another is imaginary. Furthermore, $\tilde x(s)(\dot x(s)-\mu(s,x(s),u(s,m)))$ is the only imaginary part as $ib_j\ra\tilde x(s)$. In our proposed integral we do not have any imaginary part. So we take the absolute value of it to get the magnitude of this complex number. Therefore, 
	\begin{align}\label{a10}
	\tilde x(s)(\dot x(s)-\mu(s,x(s),u(s)))&=\sqrt{\tilde x(s)^2(\dot x(s)-\mu(s,x(s),u(s)))^2}\notag\\&=\sqrt{\tilde x(s)^2}(\dot x(s)-\mu(s,x(s),u(s)))\notag\\&=|\tilde x(s)(\dot x(s)-\mu(s,x(s),u(s)))|
	\end{align}
	 Condition (\ref{a10}) implies
	\begin{align}\label{a11}
	S(x)=-\int_0^t\ \left[\sqrt{\tilde x(s)^2}(\dot x(s)-\mu(s,x(s),u(s)))-\frac{1}{2}\tilde x(s)^2 \sigma(s,x(s),u(s))\right]\ ds
	\end{align}
	The action function defined in Equation (\ref{a11}) is the integration of the Lagrangian which includes the dynamics of the system. This function and the function defined as $$ S(x)= -\ \int_0^{t}\ \left[\pi(s,x(s),u(s))+\lambda\ g(s,x(s),\dot x(s),u(s))\right]\ ds$$ are the same. Throughout this paper our one of the main objectives is to find a solution of a dynamic payoff function with an SDE as the penalization function. In this context the stochastic part appears only from the penalization function and $\lambda$ is a parameter.
	\end{proof}

\begin{prop}
	The Gaussian integral value of $\int_{\mathbb{R}}\ \exp \left\{-\frac{q}{\epsilon(1+\be)^t}\xi^2+\frac{\lambda\epsilon}{(1+\be)^t}\xi\right\}\ d\xi$ is\\ $\exp\left\{\frac{\lambda^2\epsilon^3}{4q(1+\be)^t}\right\}\sqrt{\frac{\epsilon\pi(1+\be)^t}{q}}$.
	\label{prop12}
\end{prop}

\begin{proof}
	Define $\a_1:=2q\epsilon^{-1}(1+\be)^{-t}$ and $\a_2:=\lambda\epsilon(1+\be)^{-t}$. Now we are interested in $$\int_{\mathbb{R}}\ \exp \left\{-\frac{1}{2}\a_1\xi^2+\a_2\xi\right\}d\xi.$$ 
	After factoring out yields,
	\begin{align}\label{a12}
	-\frac{1}{2}\a_1\xi^2+\a_2\xi&=-\frac{1}{2}\a_1\left(\xi^2-\frac{2\a_2}{\a_1}\xi+\frac{\a_2^2}{\a_1^2}-\frac{\a_2^2}{\a_1^2}\right)\notag\\&=-\frac{1}{2}\a_1\left(\xi-\frac{\a_2}{\a_1}\right)^2+\frac{\a_2^2}{2\a_1}
	\end{align}	
	After plugging in the result obtained in (\ref{a12}) into the integral yields
	\begin{align}\label{a13}
	\int_{\mathbb{R}}\ \exp \left\{-\frac{1}{2}\a_1\xi^2+\a_2\xi\right\}d\xi&=\int_{\mathbb{R}}\exp\left\{\frac{\a_2^2}{2\a_1}\right\}\ \exp\left\{-\frac{1}{2}\a_1\left(\xi-\frac{\a_2}{\a_1}\right)^2\right\}\ d\xi\notag\\&=\exp\left\{\frac{\a_2^2}{2\a_1}\right\}\int_{\mathbb{R}}\exp\left\{-\frac{1}{2}\a_1\left(\xi-\frac{\a_2}{\a_1}\right)^2\right\}\ d\xi\notag\\&=\exp\left\{\frac{\a_2^2}{2\a_1}\right\}\int_{\mathbb{R}}\exp\left\{-\frac{1}{2}\a_1\ \xi^2\right\}\ d\xi\notag\\&=\exp\left\{\frac{\a_2^2}{2\a_1}\right\}\ \sqrt{\frac{2\pi}{\a_1}}
	\end{align}
	After using the values of $\a_1$ and $\a_2$ we get,
	\begin{equation}\label{a14}
	\int_{\mathbb{R}}\ \exp \left\{-\frac{q}{\epsilon(1+\be)^t}\xi^2+\frac{\lambda\epsilon}{(1+\be)^t}\xi\right\}\ d\xi=\exp\left\{\frac{\lambda^2\epsilon^3}{4q(1+\be)^t}\right\}\sqrt{\frac{\epsilon\pi(1+\be)^t}{q}}.
	\end{equation}
\end{proof}

Let us consider the BPPSDE \eqref{0}.

\begin{lem}\label{ito}
	Consider the stochastic differential equation,
	\begin{align}\label{a15}
	d x &=\mu\{s,x(s),u(s),\sigma_2 [s,x(s),u(s)]\}\ ds + \{\sigma_1 (s,x(s),u(s))-\sigma_2 [s,x(s),u(s)]\}\ dW_s,
	\end{align}	
	where $\mu,u,\sigma_1$ and $\sigma_2$ are real valued functions, and $\phi^*[s,x(s),u(s)]:[0,t]\times \mathbb{R}\times[0,1]\ra\mathbb{R}$ is a continuous and at least twice differentiable in $x$ and $u$ function that is at least one time differentiable with respect to $s$ . Then
	\begin{align}
	d\phi^*[s,x(s),u(s)]&=\biggr\{\frac{\partial \phi^*[s,x(s),u(s)]}{\partial s}+\frac{\partial \phi^*[s,x(s),u(s)]}{\partial x}\mu\left\{s,x(s),u(s),\sigma_2 [s,x(s),u(s)]\right\}\notag\\&\hspace{.5cm}+\frac{1}{2}\frac{\partial^2 \phi^*[s,x(s),u(s)]}{\partial x^2}\left\{\sigma_1 [s,x(s),u(s)]-\sigma_2 [s,x(s),u(s)]\right\}^2\biggr\}ds\notag\\& \hspace{1cm}+\frac{\partial \phi^*[s,x(s),u(s)]}{\partial x}\sigma_1 [s,x(s),u(s)]dW_s\notag\\&\hspace{1.5cm}-\frac{\partial \phi^*[s,x(s),u(s)]}{\partial x}\sigma_2 [s,x(s),u(s)]\ dW_s.\notag
	\end{align}
	\label{lem13}
\end{lem}

\begin{proof}
	The general It\^o formula ($\emptyset$ksendal 2003, Theorem 4.2.1)\citep{oksendal2013stochastic} implies
	\begin{align}\label{a16}
	& d\phi^*[s,x(s),u(s)]\notag\\&=\frac{\partial \phi^*[s,x(s),u(s)]}{\partial s} ds+\frac{\partial \phi^*[s,x(s),u(s)]}{\partial x} dx(s)+\frac{\partial \phi^*[s,x(s),u(s)]}{\partial u}du(s)\notag\\&\hspace{.5cm}+\frac{1}{2}\frac{\partial^2 \phi^*[s,x(s),u(s)]}{\partial x^2}\ dx^2(s)+\frac{1}{2}\frac{\partial^2 \phi^*[s,x(s),u(s)]}{\partial x \partial u}dx(s)du(s)\notag\\&\hspace{.75cm}+\frac{1}{2}\frac{\partial^2 \phi^*[s,x(s),u(s)]}{\partial x \partial u}dx(s)du(s)+\frac{1}{2}\frac{\partial^2 \phi^*[s,x(s),u(s)]}{\partial u^2}du^2(s)
	\end{align}
	Let us assume the stubbornness is somewhat neutral over time such that $du(s)=0$. The \emph{stubbornness neutrality} means a player is consistent with their stubbornness, if there is any inconsistency then, they would not be able to make through the team. Therefore equation (\ref{a16})	becomes,
	\begin{align}\label{a17}
	d\phi^*[s,x(s),u(s)]&=\frac{\partial \phi^*[s,x(s),u(s)]}{\partial s}ds+\frac{\partial \phi^*[s,x(s),u(s)]}{\partial x}d x(s)+\frac{1}{2}\frac{\partial^2 \phi^*[s,x(s),u(s)]}{\partial x^2}\ dx^2(s)
	\end{align}
	Substituting stochastic differential equation given in Equation (\ref{a15}) into Equation (\ref{a16}) yields,
	\begin{align}\label{a18}
	& d\phi^*[s,x(s),u(s)]\notag\\&=\frac{\partial \phi^*[s,x(s),u(s)]}{\partial s}ds\notag\\&\hspace{.5cm}+\frac{\partial \phi^*[s,x(s),u(s)]}{\partial x}\biggr\{\mu\left\{s,x(s),u(s),\sigma_2 [s,x(s),u(s)]\right\}\ ds\notag\\&\hspace{1cm} + \left\{\sigma_1 [s,x(s),u(s)]-\sigma_2 [s,x(s),u(s)]\right\}\ dW_s\biggr\}\notag\\&\hspace{1.5cm}+\frac{1}{2}\frac{\partial^2 \phi^*[s,x(s),u(s)]}{\partial x^2}\biggr\{\mu\left\{s,x(s),u(s),\sigma_2 [s,x(s),u(s)]\right\}\ ds \notag\\&\hspace{2cm}+ \left\{\sigma_1 [s,x(s),u(s)]-\sigma_2 [s,x(s),u(s)]\right\}\ dW_s\biggr\}^2.
	\end{align}
	Therefore,
	\begin{align}\label{a19}
	& d\phi^*[s,x(s),u(s)]\notag\\&=\frac{\partial \phi^*[s,x(s),u(s)]}{\partial s}ds+\frac{\partial \phi^*[s,x(s),u(s)]}{\partial x}\mu\left\{s,x(s),u(s),\sigma_2 [s,x(s),u(s)]\right\}\ ds\notag\\&\hspace{.25cm}+\frac{\partial \phi^*[s,x(s),u(s)]}{\partial x} \left\{\sigma_1 [s,x(s),u(s)]-\sigma_2 [s,x(s),u(s)]\right\}\ dW_s\notag\\&\hspace{.5cm}+\frac{1}{2}\frac{\partial^2 \phi^*[s,x(s),u(s)]}{\partial x^2}\ \mu^2\left\{s,x(s),u(s),\sigma_2 [s,x(s),u(s)]\right\}\ ds^2 \notag\\&\hspace{.75cm}+\frac{1}{2}\frac{\partial^2 \phi^*[s,x(s),u(s)]}{\partial x^2}\left\{\sigma_1 [s,x(s),u(s)]-\sigma_2 [s,x(s),u(s)]\right\}^2\ dW_s^2\notag\\&\hspace{1cm}+\frac{\partial^2 \phi^*[s,x(s),u(s)]}{\partial x^2}\mu\left\{s,x(s),u(s),\sigma_2 [s,x(s),u(s)]\right\}\notag\\&\hspace{1.25cm}\left\{\sigma_1 [s,x(s),u(s)]-\sigma_2 [s,x(s),u(s)]\right\}\ ds\ dW_s.
	\end{align}
	The differential rules from It\^o's formula imply $ds^2=ds$, $ dW_s=0$ and $dW_s^2=ds$. Equation (\ref{a19}) yields,
	\begin{align}\label{a20}
	& d\phi^*[s,x(s),u(s)]\notag\\&=\frac{\partial \phi^*[s,x(s),u(s)]}{\partial s}ds+\frac{\partial \phi^*[s,x(s),u(s)]}{\partial x}\ \mu\left\{s,x(s),u(s),\sigma_2 [s,x(s),u(s)]\right\}\ ds\notag\\&\hspace{1cm}+\frac{1}{2}\frac{\partial^2 \phi^*[s,x(s),u(s)]}{\partial x^2}\left\{\sigma_1 [s,x(s),u(s)]-\sigma_2 [s,x(s),u(s)]\right\}^2\ ds\notag\\&\hspace{1.5cm}+\frac{\partial \phi^*[s,x(s),u(s)]}{\partial x} \left\{\sigma_1 [s,x(s),u(s)]-\sigma_2 [s,x(s),u(s)]\right\}\ dW_s,
	\end{align}
	giving the result. 
\end{proof}	

\begin{prop}
	(Burkholder-Davis-Gundy inequality) Let $\{\mathcal{M}\}_s=\{\phi^*(s,x(s),u(s))
	,\\ {\sigma_2}(s,x(s),u(s))\}$ be a continuous local martingale such that $\mathcal{M}_0=0$ and $s\in[0,t]$. There exists two constants $\zeta_\rho$ and $\zeta$ such that,
	\begin{align}
	& \zeta_\rho\biggr[\E\ \int_0^{\aleph_k}\ ||\phi^*(s,x(s),u(s))||_{\mathcal{H}_2}^2\ ||\hat\sigma(s,x(s),u(s))||_{\mathcal{H}_2^*}^2 \biggr]^{\frac{\delta}{2}}\notag\\&\hspace{.5cm}\leq\E\ \biggr|\sup_{s\leq\aleph_k}\ \int_0^t\ (\phi^*(s,x(s),u(s))
	,{\sigma_2}(s,x(s),u(s)))\ dW_s\biggr|^\delta\notag\\&\hspace{1cm}\leq\zeta\biggr[\E\ \int_0^{\aleph_k}\ ||\phi^*(s,x(s),u(s))||_{\mathcal{H}_2}^2\ ||\hat\sigma(s,x(s),u(s))||_{\mathcal{H}_2^*}^2 \biggr]^{\frac{\delta}{2}}
	\end{align}
	where $\delta\in[0,\infty)$ and $\zeta_\rho$ and $\zeta$ are independent of $t>0$ and $\{\mathcal{M}\}_{s\in[0,t]}$. Furthermore, when $\delta=1$ then we get ,
	\begin{align}
	& \zeta_\rho\biggr[\E\ \int_0^{\aleph_k}\ ||\phi^*(s,x(s),u(s))||_{\mathcal{H}_2}^2\ ||\hat\sigma(s,x(s),u(s))||_{\mathcal{H}_2^*}^2 \biggr]^{\frac{1}{2}}\notag\\&\hspace{.5cm}\leq\E\ \biggr|\sup_{s\leq\aleph_k}\ \int_0^t\ (\phi^*(s,x(s),u(s))
	,{\sigma_2}(s,x(s),u(s)))\ dW_s\biggr|\notag\\&\hspace{1cm}\leq\zeta\biggr[\E\ \int_0^{\aleph_k}\ ||\phi^*(s,x(s),u(s))||_{\mathcal{H}_2}^2\ ||\hat\sigma(s,x(s),u(s))||_{\mathcal{H}_2^*}^2 \biggr]^{\frac{1}{2}}\notag
	\end{align}
	\label{prop14}
\end{prop}

\begin{proof}
	Suppose, $\{\mathcal{M}\}_{s\in[0,t]}=\{\phi^*(s,x(s),u(s))
	, {\sigma_2}(s,x(s),u(s))\}$ be a continuous local martingale and $\langle\mathcal{M}\rangle_{s\in[0,t]}=||\phi^*(s,x(s),u(s))||_{\mathcal{H}_2}^2\ ||\hat\sigma(s,x(s),u(s))||_{\mathcal{H}_2^*}^2$. Let us define $\delta_1>2$ and $\kappa\in(0,1)$. It\^o formula \cite{oksendal2013stochastic} on $\{\mathcal{M}\}_{s\in[0,t]}$ implies
	\begin{align}\label{a22}
	d\ |\mathcal{M}_{s\in[0,t]}|^{\delta_1} &= \delta_1\ |\mathcal{M}_{s\in[0,t]}|^{\delta_1-1}\ \sgn (\mathcal{M}_{s\in[0,t]})\ d\mathcal{M}_{s\in[0,t]}+ \frac{1}{2}\delta_1(\delta_1-1) |\mathcal{M}_{s\in[0,t]}|^{\delta_1-2}\ d\langle \mathcal{M}_{s\in[0,t]}\rangle\notag\\&=\delta_1 \ \sgn(\mathcal{M}_{s\in[0,t]})|\ \mathcal{M}_{s\in[0,t]}|^{\delta_1-1}\ d\mathcal{M}_{s\in[0,t]}+ \frac{1}{2}\delta_1(\delta_1-1) |\mathcal{M}_{s\in[0,t]}|^{\delta_1-2}\ d\langle \mathcal{M}_{s\in[0,t]}\rangle.
	\end{align}
	 For every bounded stopping time $\upsilon$ Doob's stopping theorem yields
	\begin{align}\label{a23}
	\E\ \left\{|\mathcal{M}_{\upsilon\in[0,t]}|^{\delta_1}|\mathcal{F}_0\right\}\leq \ \frac{1}{2}\delta_1(\delta_1-1)\ \E\ \biggr\{\int_0^\upsilon\ |\mathcal{M}_{s\in[0,t]}|^{\delta_1-2}\ d\langle \mathcal{M}_{s\in[0,t]}\rangle \bigg|\mathcal{F}_0\biggr\}.
	\end{align}
	Finally, Lenglart's domination inequality implies
	\begin{align}\label{a24}
	\E\ \biggr[\biggr(\sup_{s\in(0,t)}\ |\mathcal{M}_{s\in[0,t]}|^{\delta_1}\biggr)^\kappa \biggr]\leq\frac{2-\kappa}{1-\kappa} \biggr(\frac{1}{2}\delta_1(\delta_1-1)\biggr)^\kappa\ \E\ \biggr[ \biggr(\int_0^t\ |\mathcal{M}_{s\in[0,t]}|^{\delta_1-2}\ d\langle \mathcal{M}_{s\in[0,t]}\rangle |\mathcal{F}_0\biggr)^\kappa \biggr].
	\end{align}
	Now in inequality (\ref{a24}) if we find out the bound of the expectation part, then we can get the upper bound of this entire inequality. Therefore,
	\begin{align}\label{a25}
	&\E\ \biggr[ \biggr(\int_0^t\ |\mathcal{M}_{s\in[0,t]}|^{\delta_1-2}\ d\langle \mathcal{M}_{s\in[0,t]}\rangle |\mathcal{F}_0\biggr)^\kappa \biggr]\notag\\&\hspace{1cm}\leq \E\ \biggr[\biggr(\sup_{s\in[0,t]} |\mathcal{M}_{s\in[0,t]}|\biggr)^{\kappa(\delta_1-2)}  \biggr(\int_0^t\  d\langle \mathcal{M}_{s\in[0,t]}\rangle |\mathcal{F}_0\biggr)^\kappa \biggr]\notag\\&\hspace{2cm}\leq \E\ \biggr[\biggr(\sup_{s\in[0,t]} |\mathcal{M}_{s\in[0,t]}|\biggr)^{\kappa\delta_1}\biggr]^{1-\frac{2}{\delta_1}}\ \E\ \biggr[\langle \mathcal{M}_t\rangle^{\frac{\kappa \delta_1}{2}}\biggr]^{\frac{2}{\delta_1}}.
	\end{align}
	After using the result of the inequality (\ref{a25}) into (\ref{a24}) we get the upper bound of the expectation as 
	\begin{align}\label{a26}
	&\E\ \biggr[\biggr(\sup_{s\in[0,t]}\ |\mathcal{M}_{s\in[0,t]}|^{\delta_1}\biggr)^\kappa\biggr]\notag\\&\hspace{1cm}\leq\frac{2-\kappa}{1-\kappa} \biggr(\frac{1}{2}\delta_1(\delta_1-1)\biggr)^\kappa \E\ \biggr[\biggr(\sup_{s\in[0,t]} |\mathcal{M}_{s\in[0,t]}|\biggr)^{\kappa\delta_1}\biggr]^{1-\frac{2}{\delta_1}} \E\ \biggr[\langle \mathcal{M}_t\rangle^{\frac{\kappa \delta_1}{2}}\biggr]^{\frac{2}{\delta_1}}.
	\end{align}
	Now as $\delta=\delta_1 \kappa$, for $s\leq\aleph_k\in[0,t]$ we get,
	\begin{align}\label{a27}
	&\E\ \biggr[\sup_{s\leq\aleph_k}\ |\mathcal{M}_{s\in[0,t]}|^{\delta}\biggr]\leq \zeta \E \biggr[\langle\mathcal{M}_t\rangle^\frac{\delta}{2}\biggr]\notag\\&\hspace{.5cm}\implies  \E\ \biggr|\sup_{s\leq\aleph_k}\ \int_0^t\ (\phi^*(s,x(s),u(s))
	,{\sigma_2}(s,x(s),u(s)))\ dW_s\biggr|^\delta\notag\\&\hspace{1cm}\leq\zeta\biggr[\E\ \int_0^{\aleph_k}\ ||\phi^*(s,x(s),u(s))||_{\mathcal{H}_2}^2\ ||\hat\sigma(s,x(s),u(s))||_{\mathcal{H}_2^*}^2 \biggr]^{\frac{\delta}{2}}.
	\end{align}
	Furthermore, when $\delta=1$, we have
	\begin{align}\label{a28}
	&\E\ \biggr|\sup_{s\leq\aleph_k}\ \int_0^t\ (\phi^*(s,x(s),u(s))
	,{\sigma_2}(s,x(s),u(s)))\ dW_s\biggr|\notag\\&\hspace{1cm}\leq\zeta\biggr[\E\ \int_0^{\aleph_k}\ ||\phi^*(s,x(s),u(s))||_{\mathcal{H}_2}^2\ ||\hat\sigma(s,x(s),u(s))||_{\mathcal{H}_2^*}^2 \biggr]^{\frac{1}{2}}.
	\end{align}
	Therefore, inequality (\ref{a28}) shows the right inequality of the proposition. Oksendal (2013) \cite{oksendal2013stochastic} implies
	\begin{align}\label{a29}
	\mathcal{M}_{s\in[0,t]}^2=\langle \mathcal{M}_{s\in[0,t]}\rangle +2 \int_0^t\ \mathcal{M}_{s\in[0,t]}\ d\mathcal{M}_{s\in[0,t]}.
	\end{align}
	Thus,
	\begin{align}\label{a30}
	\E\ \big[\langle \mathcal{M}_{t}\rangle^\frac{\delta}{2}\big]\leq\ \zeta_\rho\biggr( \E\ \biggr[\big(\sup_{s\in[0,t]}\ |\mathcal{M}_{s\in[0,t]}|\big)^\delta\biggr]+\ \E\ \biggr[\sup_{s\in[0,t]}\ \biggr|\int_0^t\ \mathcal{M}_{s\in[0,t]}\ d\mathcal{M}_{s\in[0,t]}\biggr|^{\frac{\delta}{2}}\biggr]\biggr).
	\end{align}
	Similarly,
	\begin{align}\label{a31}
	&\E\ \biggr[ \sup_{s\leq\aleph_k}\biggr|\int_0^s \mathcal{M}_{q\in[0,t]}\ d\mathcal{M}_{q\in[0,t]} \biggr|^{\frac{\delta}{2}}\biggr]\leq C_\delta\ \E\biggr[\biggr(\int_0^{\aleph_k}\ \mathcal{M}_{q\in[0,t]}^2 \ d\langle\mathcal{M}_{q\in[0,t]}\rangle\biggr)^{\frac{\delta}{4}}\biggr]\notag\\&\hspace{.5cm}\leq C_\delta\ \E\biggr[ \biggr(\sup_{s\leq\aleph_k}\big|\mathcal{M}_{q\in[0,t]}\big|^{\frac{\delta}{2}}\ \langle \mathcal{M}_{\aleph_k}\rangle^{\frac{\delta}{4}}\biggr)\biggr]\notag\\&\hspace{1cm}\leq  C_\delta\ \E\biggr[ \biggr(\sup_{s\leq\aleph_k}\big|\mathcal{M}_{q\in[0,t]}\big|\biggr)^\delta \biggr]^{\frac{1}{2}}\ \E\ \biggr[\langle \mathcal{M}_{\aleph_k}\rangle \biggr]^{\frac{1}{2}}.
	\end{align}
	Therefore, 
	\begin{align}\label{a32}
	\E\biggr[\langle \mathcal{M}_{\aleph_k}\rangle^{\frac{1}{2}}\biggr]&\leq \zeta_\rho\ \biggr(\E\ \biggr[ \big(\sup_{s\in[0,t]}\mathcal{M}_{s\in[0,t]}|\mathcal{M}_{s\in[0,t]}|\big)^\delta\biggr]\notag\\&\hspace{1cm}+C_\delta\ \E\biggr[ \biggr(\sup_{s\leq\aleph_k}\big|\mathcal{M}_{q\in[0,t]}\big|\biggr)^\delta \biggr]^{\frac{1}{2}}\ \E\ \biggr[\langle \mathcal{M}_{\aleph_k}\rangle \biggr]^{\frac{1}{2}}\biggr).
	\end{align}
	If we carefully look at inequality (\ref{a32}) we find out that it is in the form of $m^2\leq \zeta_\rho (n^2+C_\delta\ mn)$ which further implies $\zeta_\rho m^2\leq n^2 $ [ as $2mn\leq \frac{1}{\epsilon}\ m^2 +\epsilon n^2$], for any chosen $\epsilon$. Hence,
	\begin{align}\label{a33}
	& \zeta_\rho\biggr[\E\ \int_0^{\aleph_k}\ ||\phi^*(s,x(s),u(s))||_{\mathcal{H}_2}^2\ ||\hat\sigma(s,x(s),u(s))||_{\mathcal{H}_2^*}^2 \biggr]^{\frac{\delta}{2}}\notag\\&\hspace{.5cm}\leq\E\ \biggr|\sup_{s\leq\aleph_k}\ \int_0^t\ (\phi^*(s,x(s),u(s))
	,{\sigma_2}(s,x(s),u(s)))\ dW_s\biggr|^\delta,
	\end{align}
	with $\delta=1$ we get,
	\begin{align}\label{a34}
	& \zeta_\rho\biggr[\E\ \int_0^{\aleph_k}\ ||\phi^*(s,x(s),u(s))||_{\mathcal{H}_2}^2\ ||\hat\sigma(s,x(s),u(s))||_{\mathcal{H}_2^*}^2 \biggr]^{\frac{1}{2}}\notag\\&\hspace{.5cm}\leq\E\ \biggr|\sup_{s\leq\aleph_k}\ \int_0^t\ (\phi^*(s,x(s),u(s))
	,{\sigma_2}(s,x(s),u(s)))\ dW_s\biggr|.
	\end{align}
	Inequalities (\ref{a34}) and (\ref{a28}) imply
	\begin{align}\label{a35}
	& \zeta_\rho\biggr[\E\ \int_0^{\aleph_k}\ ||\phi^*(s,x(s),u(s))||_{\mathcal{H}_2}^2\ ||\hat\sigma(s,x(s),u(s))||_{\mathcal{H}_2^*}^2 \biggr]^{\frac{1}{2}}\notag\\&\hspace{.5cm}\leq\E\ \biggr|\sup_{s\leq\aleph_k}\ \int_0^t\ (\phi^*(s,x(s),u(s))
	,{\sigma_2}(s,x(s),u(s)))\ dW_s\biggr|\notag\\&\hspace{1cm}\leq\zeta\biggr[\E\ \int_0^{\aleph_k}\ ||\phi^*(s,x(s),u(s))||_{\mathcal{H}_2}^2\ ||\hat\sigma(s,x(s),u(s))||_{\mathcal{H}_2^*}^2 \biggr]^{\frac{1}{2}}.
	\end{align}
	This completes the proof.
\end{proof}

\medskip

\begin{prop}
	(Grownwall inequality) Let us assume $\mathcal{H}_1$ be a Banach space such that there exists an open subset $\mathcal{S}_b$ such that $\mathcal{S}_b\subset\mathcal{H}_1$. Suppose, there exists two continuous functions such that $f_1,f_2:[\a,\be]\times\mathcal{S}_b\ra\mathcal{H}_1$ and $m_1,m_2:[\a,\be]\ra\mathcal{S}_b$ satisfy the initial value problems
	\begin{align}
	m_1'(s)&=f_1(s,m_1(s)),\ m_1(\a)=m_{10},\notag\\ m_2'(s)&=f_2(s,m_2(s)),\ m_2(\a)=m_{20}.\notag
	\end{align} 
	There exists a constant $\zeta$ such that,
	\begin{align}
	||f_2(s,e_2)-f_2(s,e_1)||\leq\zeta\ ||e_2-e_1||,\notag 
	\end{align}
	and a continuous function $\mho:[\a,\be]\ra[0,\infty)$ so that
	\begin{align}
	||f_1(s,m_1(s))-f_2(s,m_1(s))||\leq\mho(s).\notag
	\end{align}
	Then
	\begin{align}
	||m_1(t)-m_2(t)||\leq\ e^{\zeta|t-\a|}\ ||m_{10}-m_{20}||+e^{\zeta|t-\a|}\ \int_\a^t\ e^{-\zeta|s-\a|}\ \mho(s)\ ds,\notag 
	\end{align}
	where $s\in[\a,\be]$.
	\label{prop15}
\end{prop}

\begin{proof}
	For any $C^1$ function $f:[\a,\be]\mapsto\mathcal{H}_1$ we know $\frac{d}{ds}||f(s)||\leq||f'(s)||$. Consider $m_1(.),\ m_2(.):[\a,\be]^2\ra\mathcal{H}_1^2$. Then,
	\begin{align}\label{a36}
	\frac{d}{ds}||m_1(s)-m_2(s)||&\leq||m_1'(s)-m_2'(s)||\notag\\&=||f_1(s,m_1(s))-f_2(s,m_2(s))||\notag\\&\leq||f_1(s,m_1(s))-f_2(s,m_2(s))||+||f_2(s,m_1(s))-f_2(s,m_2(s))||\notag\\&\leq\mho(s)+\zeta||m_1(s)-m_2(s)||,
	\end{align}
	where $\mho(s):[\a,\be]\ra[0,\infty)$, and it is assumed that $||f_1(s,m_1(s))-f_2(s,m_2(s))||\leq\mho(s)$.
	After rearranging the inequality (\ref{a36}) we get,
	\begin{align}\label{a37}
	\frac{d}{ds}||m_1(s)-m_2(s)||-\zeta||m_1(s)-m_2(s)||\leq \mho(s).
	\end{align}
	After multiplying the integrating factor $e^{-\zeta\ s}$ in both sides of the inequality (\ref{a37}) we get,
	\begin{align}\label{a38}
	\frac{d}{ds}\ \biggr(e^{-\zeta s}||m_1(s)-m_2(s)||\biggr)\leq e^{-\zeta s}\ \mho(s).
	\end{align}
	Therefore, the integral form becomes,
	\begin{align}\label{a39}
	\int_\a^t\ \biggr[\frac{d}{ds}\ \biggr(e^{-\zeta s}||m_1(s)-m_2(s)||\biggr)\biggr]\ ds&\leq\int_\a^t e^{-\zeta s}\ \mho(s)ds\notag\\e^{-\zeta t}||m_1(t)-m_2(t)||-e^{-\zeta \a}||m_{10}-m_{20}||&\leq \int_\a^t e^{-\zeta s}\ \mho(s)ds.
	\end{align}
	Inequality (\ref{a39}) and the argument in the proposition are the same. This completes the proof.
\end{proof}

\medskip

\begin{cor}
	Let us assume $\mathcal{H}_1$ be a Banach space such that there exists an open subset $\mathcal{S}_b$ such that $\mathcal{S}_b\subset\mathcal{H}_1$. Suppose, there exists a continuous function such that $f_1:[\a,\be]\times\mathcal{S}_b\ra\mathcal{H}_1$ and $m_1,m_2:[\a,\be]\ra\mathcal{S}_b$ satisfy the initial value problems
	\begin{align}
	m_1'(s)&=f_1(s,m_1(s)),\ m_1(\a)=m_{10},\notag\\ m_2'(s)&=f_2(s,m_2(s)),\ m_2(\a)=m_{20}.\notag
	\end{align} 
	There exists a constant $\zeta\in[0,\infty)$ such that,
	\begin{align}
	||f_2(s,e_2)-f_2(s,e_1)||\leq\zeta\ ||e_2-e_1||.\notag 
	\end{align}
	Then 
	\begin{align}
	||m_1(t)-m_2(t)||\leq\ e^{\zeta|t-\a|}\ ||m_{10}-m_{20}||,
	\end{align}
	for all $t\in[\a,\be]$.
	\label{cor16}
\end{cor}

\begin{proof}
	In Proposition \ref{prop15}, assume that $f_1(.)=f_2(.)$. Since for any continuous function $\mho(s):[\a,\be]\mapsto[0,\infty)$, Proposition \ref{prop15} implies that $||f_1(s,m_1(s))-f_2(s,m_1(s))||\leq\ \mho(s)$. As $f_1$ and $f_2$ are the same then $\mho(s)\equiv0$ for all $s\in[\a,\be]$. Therefore, the second right hand term of the inequality 
	\begin{align}
	||m_1(t)-m_2(t)||\leq\ e^{\zeta|t-\a|}\ ||m_{10}-m_{20}||+e^{\zeta|t-\a|}\ \int_\a^t\ e^{-\zeta|s-\a|}\ \mho(s)\ ds\notag 
	\end{align}
	vanishes and we remain with,
	\begin{align}\label{a40}
	||m_1(t)-m_2(t)||\leq\ e^{\zeta|t-\a|}\ ||m_{10}-m_{20}||.
	\end{align}
	This completes the proof.
\end{proof}

\bibliographystyle{apalike}
\bibliography{bib}
\end{document}